\pdfoutput=1
\RequirePackage{silence}
\documentclass[a4paper,svgnames]{amsart}
\usepackage[hmarginratio={1:1},vmarginratio={1:1},lmargin=60.0pt,tmargin=59.0pt]{geometry}

\usepackage{booktabs}
\allowdisplaybreaks

\usepackage[T1]{fontenc}
\usepackage{latexsym,exscale,amsfonts,amssymb,mathtools}
\usepackage{amsmath,amsthm,amsfonts,amssymb,amscd,textcomp,bbm}
\usepackage{mathrsfs,stackrel}
\usepackage{etoolbox}
\usepackage{tikz,tikz-cd}
\usetikzlibrary{matrix,arrows.meta,decorations.markings,shapes.multipart,decorations.pathreplacing,decorations.pathmorphing,shapes.geometric}
\usepackage{aliascnt}
\usepackage{ytableau}
\usepackage{xparse}
\usepackage{cite}

\usepackage{dynkin-diagrams}
\usepackage{nicematrix}
\usepackage{tikz-3dplot}
\usepackage{mathrsfs}
\DeclareMathAlphabet{\mathscrbf}{OMS}{mdugm}{b}{n}
\usepackage{float}
\usepackage{caption}
\usepackage[most]{tcolorbox}

\usepackage{array}
\newcolumntype{C}{>{$}c<{$}}

\usepackage{xcolor,colortbl}

\definecolor{mygray}{gray}{0.6}
\definecolor{mygraydark}{gray}{0.4}
\definecolor{mygraylight}{gray}{0.85}
\definecolor{spinach}{RGB}{46,139,87}
\definecolor{tomato}{RGB}{255,99,71}
\definecolor{orchid}{RGB}{143,40,194}
\definecolor{neon}{RGB}{77,77,255}
\definecolor{lightneon}{RGB}{110,110,255}
\definecolor{pumpkin}{RGB}{224,180,80}
\definecolor{citron}{RGB}{190,180,90}

\definecolor{lava}{RGB}{207,16,32}
\definecolor{cream}{RGB}{255,253,208}
\definecolor{verdigris}{RGB}{67,179,174}
\definecolor{Black}{RGB}{0,0,0}
\definecolor{mydarkblue}{RGB}{10,10,170}
\definecolor{darkspinach}{RGB}{20,70,20}
\definecolor{darktomato}{RGB}{155,40,30}
\definecolor{darkorchid}{RGB}{50,10,100}
\definecolor{darklava}{RGB}{150,8,16}

\definecolor{zero}{RGB}{0,0,0}
\definecolor{one}{RGB}{255,0,0}
\definecolor{two}{RGB}{0,255,0}
\definecolor{three}{RGB}{0,0,255}

\usepackage{todonotes}

\usepackage{enumitem}
\setlist[enumerate]{itemsep=0.15cm,label=\emph{\upshape(\alph*)}}
\setlist[enumerate,2]{itemsep=0.15cm,label=\emph{\upshape(\roman*)}}
\setlist[enumerate,3]{itemsep=0.15cm,label=\emph{\upshape(\Alph*)}}

\let\emph\relax
\DeclareTextFontCommand{\emph}{\bfseries\em}

\renewcommand{\dots}{\text{...}}
\renewcommand{\vdots}{\rotatebox{90}{\text{...}}}
\renewcommand{\ddots}{\raisebox{0.175cm}{\rotatebox{-45}{\text{...}}}}

\newcommand{\placeholder}{{}_{-}}
\newcommand{\mystrut}{\rule[-0.2\baselineskip]{0pt}{0.9\baselineskip}}

\newcommand{\ie}{\text{i.e.}}
\newcommand{\eg}{\text{e.g.}}
\newcommand{\cf}{\text{cf.}}
\newcommand{\etc}{\text{etc.}}
\newcommand{\aka}{\text{a.k.a.}}

\usepackage[all]{xy}
\usepackage{tikz}
\usetikzlibrary{cd}
\usetikzlibrary{decorations}
\usetikzlibrary{decorations.markings}
\usetikzlibrary{decorations.pathreplacing}
\usetikzlibrary{decorations.pathmorphing}
\usetikzlibrary{arrows.meta,shapes,positioning,matrix,calc}
\usetikzlibrary{shapes.callouts}
\tikzset{anchorbase/.style={baseline={([yshift=-0.5ex]current bounding box.center)}},
tinynodes/.style={font=\tiny,text height=0.25ex,text depth=0.05ex},
smallnodes/.style={font=\scriptsize,text height=0.75ex,text depth=0.15ex},
usual/.style={line width=2.0,color=black},
}

\newcommand{\vpar}{\changed{\varstuff{q}}}

\newcommand{\C}{\mathbb{C}}
\newcommand{\R}{\mathbb{R}}
\newcommand{\Pp}{\mathbb{P}}

\newcommand{\N}{\mathbb{Z}_{\geq 0}}

\newcommand{\Z}{\mathbb{Z}}

\newcommand{\setstuff}[1]{\mathrm{#1}}

\newcommand{\varstuff}[1]{\mathtt{#1}}

\newcommand{\morstuff}[1]{\mathrm{#1}}

\newcommand{\idmor}{\morstuff{id}}

\newcommand{\sank}{n}
\newcommand{\mcoeff}{d}
\newcommand{\sym}[1][\sank]{\setstuff{S}_{#1}}

\newcommand{\br}{\leq_{B}}
\newcommand{\bbr}{<_{B}}

\newcommand{\rbb}{>_{B}}

\newcommand{\ol}{\leq_{1}}

\newcommand{\loo}{>_{1}}

\newcommand{\ldes}{\mathrm{ldes}}
\newcommand{\rdes}{\mathrm{rdes}}
\newcommand{\ep}[1][\sank]{\mathrm{EP}(#1)}

\newcommand{\kl}[1][w]{\setstuff{P}_{#1}}
\newcommand{\KL}[2]{\setstuff{P}_{#1,#2}}

\newcommand{\ssum}{{\textstyle\sum}}
\newcommand{\sspan}{\changed{\setstuff{sp}\,}}

\newcommand{\drawstringdiagram}[2]{
\begin{tikzpicture}[anchorbase]
\def\permutation{#1} % Input permutation ({\eg}, {8, 9, 10, 4, 5, 6, 7, 1, 2, 3})
\def\n{#2} % Length of the permutation
\foreach \i in {1,...,\n} {
\draw[fill=black] (\i, 0) circle (0pt); % Bottom row points
\draw[fill=black] (\permutation[\i-1], 2) circle (0pt); % Top row points
}

\foreach \i in {1,...,\n} {
\draw[usual] (\i, 0) -- (\permutation[\i-1], 2);
}
\end{tikzpicture}
}

\newcommand{\drawstringdiagramlabel}[2]{
\begin{tikzpicture}[anchorbase]
\def\permutation{#1} % Input permutation (e.g., [8, 9, 10, 4, 5, 6, 7, 1, 2, 3])
\def\n{#2} % Length of the permutation
\foreach \i in {1,...,\n} {
    \draw[fill=black] (\i, 0) circle (0pt); % Bottom row points
    \draw[fill=black] (\i, 2) circle (0pt); % Top row points (now at positions 1 to n)
}

\foreach \i in {1,...,\n} {
    \draw[usual] (\i, 0) -- (\permutation[\i-1], 2);
    \node at (\i,-0.5) {\i};
    \node at (\i,2.5) {\i};
}
\end{tikzpicture}
}

\def\NewTheorem#1{%
\newaliascnt{#1}{equation}%
\newtheorem{#1}[#1]{#1}%
\aliascntresetthe{#1}%
\expandafter\def\csname #1autorefname\endcsname{#1}%
}
\def\equationautorefname~#1\null{(#1)\null}

\numberwithin{equation}{subsection}

\NewTheorem{Proposition}
\NewTheorem{Theorem}
\NewTheorem{Corollary}
\AtEndEnvironment{Corollary}{\null\hfill$\square$}%
\NewTheorem{Lemma}
\NewTheorem{Conjecture}
\NewTheorem{Speculation}
\NewTheorem{Question}
\NewTheorem{Assumption}
\theoremstyle{definition}
\NewTheorem{Definition}
\AtEndEnvironment{Definition}{\null\hfill$\Diamond$}%
\NewTheorem{Classification Problem}
\AtEndEnvironment{Classification Problem}{\null\hfill$\Diamond$}%
\NewTheorem{Notation}
\AtEndEnvironment{Notation}{\null\hfill$\Diamond$}%
\NewTheorem{Example}
\AtEndEnvironment{Example}{\null\hfill$\Diamond$}%
\NewTheorem{Examples}
\AtEndEnvironment{Examples}{\vskip-10mm\null\hfill$\Diamond$}%

\theoremstyle{remark}
\NewTheorem{Remark}
\AtEndEnvironment{Remark}{\null\hfill$\Diamond$}%

\usepackage[hypertexnames=false]{hyperref}
\usepackage{bookmark}
\hypersetup{
pdftoolbar=true,
pdfmenubar=true,
pdffitwindow=false,
pdfstartview={FitH},
pdftitle={Big data approach to Kazhdan--Lusztig polynomials},
pdfauthor={Abel Lacabanne, Daniel Tubbenhauer and Pedro Vaz},
pdfsubject={},
pdfcreator={Abel Lacabanne, Daniel Tubbenhauer and Pedro Vaz},
pdfproducer={Abel Lacabanne, Daniel Tubbenhauer and Pedro Vaz},
pdfkeywords={},
pdfnewwindow=true,
colorlinks=true,
linkcolor=mydarkblue,
citecolor=teal,
filecolor=magenta,
urlcolor=orchid,
linkbordercolor=lava,
citebordercolor=teal,
urlbordercolor=orchid,
linktocpage=true
}

\def\makeautorefname#1#2{\csdef{#1autorefname}{#2}}
\makeautorefname{section}{Section}%
\makeautorefname{subsection}{Section}%
\makeautorefname{subsubsection}{Section}%

\def\changed#1{{#1}}
\def\ochanged#1{{#1}}

\begin{document}
\title[Big data approach to Kazhdan--Lusztig polynomials]{Big data approach to Kazhdan--Lusztig polynomials}
\author[A. Lacabanne, D. Tubbenhauer and P. Vaz]{Abel Lacabanne, Daniel Tubbenhauer and Pedro Vaz}

\address{A.L.: Laboratoire de Math{\'e}matiques Blaise Pascal (UMR 6620), Universit{\'e} Clermont Auvergne, Campus Universitaire des C{\'e}zeaux, 3 place Vasarely, 63178 Aubi{\`e}re Cedex, France, \href{http://www.normalesup.org/~lacabanne}{www.normalesup.org/$\sim$lacabanne},
\href{https://orcid.org/0000-0001-8691-3270}{ORCID 0000-0001-8691-3270}}
\email{abel.lacabanne@uca.fr}

\address{D.T.: The University of Sydney, School of Mathematics and Statistics F07, Office Carslaw 827, NSW 2006, Australia, \href{http://www.dtubbenhauer.com}{www.dtubbenhauer.com}, \href{https://orcid.org/0000-0001-7265-5047}{ORCID 0000-0001-7265-5047}}
\email{daniel.tubbenhauer@sydney.edu.au}

\address{P.V.: Institut de Recherche en Math{\'e}matique et Physique, 
Universit{\'e} catholique de Louvain, Chemin du Cyclotron 2,  
1348 Louvain-la-Neuve, Belgium, \href{https://perso.uclouvain.be/pedro.vaz}{https://perso.uclouvain.be/pedro.vaz}, \href{https://orcid.org/0000-0001-9422-4707}{ORCID 0000-0001-9422-4707}}
\email{pedro.vaz@uclouvain.be}

\begin{abstract}
We investigate the structure of Kazhdan--Lusztig polynomials of the symmetric group by leveraging computational approaches from big data, including exploratory and topological data analysis, applied to the polynomials for symmetric groups of up to 11 strands.
\end{abstract}


\addtocontents{toc}{\protect\setcounter{tocdepth}{1}}

\maketitle

\tableofcontents

\section{Introduction}\label{S:Intro}

We investigate a data set of Kazhdan--Lusztig (KL) polynomials using techniques traditionally employed in data science.

\subsection{Background and ideas}

The \emph{Kazhdan–Lusztig (KL) polynomials} are fundamental yet enigmatic objects in combinatorial representation theory. First introduced in \cite{KaLu-reps-coxeter-groups} for arbitrary Coxeter groups, this paper focuses on the type A Coxeter group, the symmetric group $\sym$ on $\{1,\dots,\sank\}$. In this context, KL polynomials represent entries in a nonnegatively graded change-of-basis matrix between simple and Verma modules of $\mathfrak{sl}_{\sank}$. As a result, these polynomials are either zero or belong to $1+\vpar\N[\vpar]$, where $\vpar$ denotes the grading variable.

Although a main focus of research in combinatorics, geometry, and representation theory alike, \changed{these polynomials still hold many mysteries.}
The starting point of this work is the observation that KL polynomials exhibit patterns akin to statistical distributions. By analyzing these distributions, we aim to uncover structures that remain elusive through traditional methods, such as combinatorial or geometric approaches.
Our methodology involves systematic data processing and analysis, often referred to as \emph{``big data,''} utilizing techniques such as data visualization, exploratory data analysis (EDA), and topological data analysis (TDA).
Our approach is inspired by work in knot theory as, for example, \cite{LeHaSa-jones,DlGuSa-mapper,TuZh} 
and differs from the deep learning approaches to KL theory as in {\eg} \cite{Nature,Wi-deep}.

Using these methods, we will discuss several conjectures about KL polynomials and their distribution, \changed{and for some of them, we provide a possible approach for proving these conjectures.} \ochanged{Additionally, the data science approach offers a different perspective, revealing patterns that, while possibly beyond formal verification, are still worth highlighting.}

\begin{Remark}
A critical aspect of this study is the scale of the data: KL polynomials are indexed by pairs of elements in $\sym$, resulting in $(\sank!)^2$ \ochanged{values}. Although computable, their computation for $\sym[11]$, with one permutation fixed as trivial, required approximately 60 days using the program in \cite{Wa-mu-for-sn} on the servers of Laboratoire de Math{\'e}matiques Blaise Pascal, Universit{\'e} Clermont Auvergne. For $\sym[12]$, we anticipate significantly higher computational costs, at least 12 times longer, although this is a very conservative lower bound that completely excludes potential memory limitations. \changed{One particularly notable KL polynomial, with potentially the largest value at $\vpar=1$ (when normalized by dividing by the group size) ever computed, see \autoref{E:Big} below, is for $\sym[13]$ and required approximately 60 days of computation on the same server.}
\end{Remark}

\begin{Notation}
In this paper, we have \emph{conjectures} and \emph{speculations}. 
Conjectures are presented in their standard sense, while speculations refer to preliminary hypotheses that lack full support from the data.
We hope that both serve as an inspiration to prove or, {\color{purple}equally exciting}, disprove the corresponding statements\ochanged{, or simply to spark new ideas}. We also have \emph{questions} meant in the standard way of the word.
\end{Notation}

\begin{Remark}\label{R:Code}
All the data files and supplementary material 
(such as higher resolution pictures, code, a possible empty Erratum {\etc}) can be found online at \cite{LaTuVa-kl-big-data-code}.
\end{Remark}

\subsection{Main results and their philosophy}

\changed{Consider the following related yet distinct questions:}
\begin{enumerate}[label=(\roman*)]

\item \changed{What is the distribution of algebraic integers in $\C$ (identified with the real plane)?}

\item \changed{For a fixed $\sank \in \N$, what is the distribution of algebraic integers in $\C$ arising from polynomials of degree at most $\sank$?}

\end{enumerate}
\changed{Both questions have a different focus and are legitimate, and, surprisingly, their answers diverge sharply. The first is straightforward: since any plot has finitely many pixels, rendering algebraic integers in black against a white background produces a fully black image. This happens since they are dense in $\C$. The second, however, reveals a striking phenomenon; fractal patterns emerge, as exemplified on the first page of \cite{BaChDe-roots}. The distribution is highly nonuniform and invites deeper study, though many observed patterns lack rigorous proof, relying instead on experimental mathematics to uncover their properties.}

\begin{figure}[H]
\fcolorbox{tomato!50}{orchid!10}{\includegraphics[height=6.0cm]{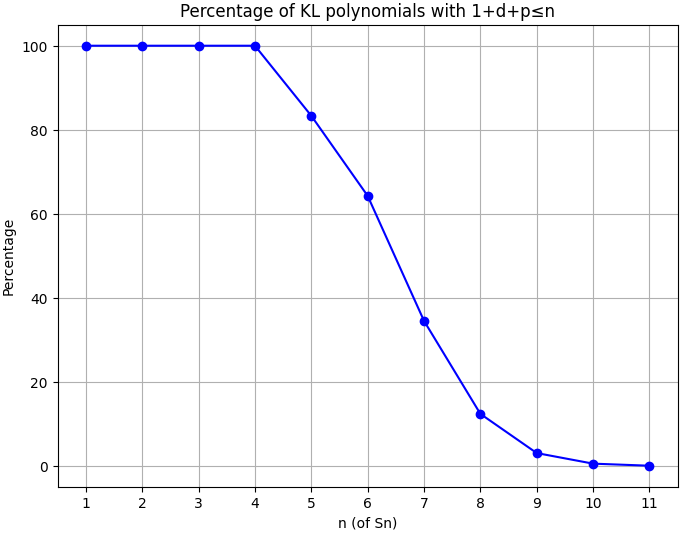}\includegraphics[height=6.0cm]{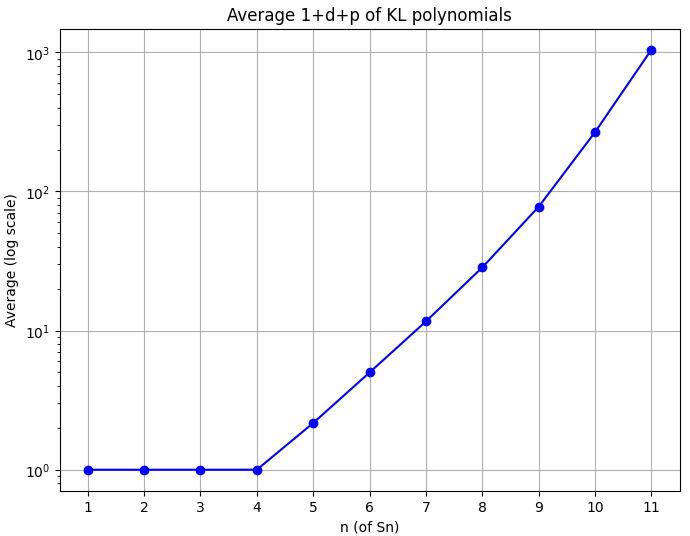}}
\begin{gather*}
\begin{tabular}{c|ccccccccccc}
$\sank$  & 1 & 2 & 3 & 4 & 5 & 6 & 7 & 8 & 9 & 10 & 11 \\
\hline
\% & 100 & 100 & 100 & 83.33 & 64.28 & 64.29 & 34.54 & 12.41 & 3.09 & 0.57 & 0.08 \\
av. $1+d+p$ & 1 & 1 & 1 & 1.3333 & 2.1666 & 5 & 11.6363 & 28.4379 & 77.6 & 266.6218 & 1043.3541 \\
\end{tabular}
\end{gather*}
\caption{The percentage of KL polynomials $\kl$ (with one trivial, as defined in \autoref{S:Notation}) for $w\in\sym$ satisfying $1+d+p\leq\sank$, and the average $1+d+p$, in log scale, running over all such KL polynomials. Here $d$ is the degree and $p$ is the evaluation at $\vpar=1$.}
\label{polo}
\end{figure}

\changed{For KL polynomials, \cite{Po-kl-symmetric-group} established a foundational and beautiful result: every polynomial in $1+\vpar\N[\vpar]$ arises as a KL polynomial for some $\sym$ (``every polynomial is a KL polynomial'') through an algorithm constructing such polynomials. This mirrors the broad scope of (i). By contrast, this paper explores the analog of (ii), examining KL polynomials under rank constraints. This departs significantly from \cite{Po-kl-symmetric-group}: its algorithm (the only known general construction, to our knowledge) generates a polynomial of degree $d$ and evaluation $p$ at $\vpar=1$ within $\sym[{1+d+p}]$. This imposes a severe bias, with very few polynomials satisfying this bound; see \autoref{polo}. 
For instance, consider the polynomial in $\sym[13]$ from \autoref{E:Big}. To realize it, the algorithm in \cite{Po-kl-symmetric-group} requires $\sym[67994199]$, a group whose order, $67994199!$, exceeds 500000000 digits.
In that sense, ``every polynomial is a KL polynomial'' may be a misleading slogan. Rather, we ask:}

\begin{gather*}
\fcolorbox{orchid!50}{spinach!10}{\mystrut``What is the distribution of KL polynomials for fixed $\sank$?''} 
\end{gather*}

\changed{Here, \emph{distribution} should be interpreted in a statistical or probabilistic context. Specifically, analyzing a distribution involves, for instance, examining its support or (probability) density, its extreme values, and its overall structure. Accordingly, we divide the analysis of the distribution into three distinct components: the study of the \emph{density} of the KL polynomials in \autoref{S:Density} and \autoref{S:Number}, the \emph{extremes} in \autoref{S:Growth}, and the \emph{structural properties} in \autoref{S:Unimodal}--\autoref{S:Ballmapper}.}

\changed{Our main results, largely conjectural and driven by \emph{experimental methods rooted in data}, then include the following, alongside smaller findings.}
\begin{enumerate}[label=$\blacktriangleright$]

\item \changed{\textit{Density.} Almost all KL polynomials for $\sym$ are zero, even on average; see \autoref{S:Density}. Yet, the average decreases much more gradually than the proportion of nonzero polynomials. However, and surprising, when focusing on KL polynomials with one permutation fixed as trivial, then almost all of them are not just nonzero, but also different, {\cf} \autoref{S:Number}.}

\item \changed{\textit{Extremes.} Some KL polynomials exhibit exceptionally large coefficients. We conjecture in \autoref{S:Growth} that their growth is superexponential in the rank, far exceeding the linear construction of \cite{Po-kl-symmetric-group}. For example, at $\sank=13$, a coefficient reaches eight digits; see \autoref{E:Big}.}

\item \changed{\textit{Structure.} The distribution of KL polynomials reveals a distinctive and surprising structure: nearly all are unimodal, with roots exhibiting fractal-like patterns and other features; see \autoref{S:Unimodal}--\autoref{S:Ballmapper} for details. This potentially explains why the algorithm from \cite{Po-kl-symmetric-group} may need to be rather inefficient: a randomly selected KL polynomial differs markedly from a random polynomial, and the set of KL polynomials is so biased that certain polynomials appear very rarely and only for $\sank\gg0$.}

\end{enumerate}
\changed{Where possible, we provide evidence supporting these conjectured patterns.}

\subsection{Future directions}

The immediate candidates that come to mind to apply big data methods are graph and knot polynomials,
but this section focuses on potential representation-theoretical applications of big data methods. We outline several avenues for exploration:
\begin{enumerate}[label=(\Roman*)]

\item \textit{Generalizations to other KL polynomials.}
An immediate extension is to study variations of KL polynomials, such as other types, antispherical or spherical KL polynomials, which have deep connections to representation theory, {\cf} \cite{So-tilting-a}. Data for such cases can potentially be generated using existing software, {\eg} \cite{duCl-positivity-finite-hecke,Gi-ihecke-code}. Another intriguing direction is to investigate the \emph{p-canonical (pKL) polynomials}, {\cf} \cite{JeWi-p-canonical}, though their computational complexity presents a significant challenge. Special cases, such as those in \cite{Ba-pkl}, may provide a tractable starting point for data collection and analysis.

\item \textit{Canonical bases and change-of-basis polynomials.}
Similarly and related, {\cf} \cite{FrKhKi-kl-and-uqsl2}, \emph{canonical ({\aka} crystal) bases} of Lie algebras and their change-of-basis matrices give interesting polynomials. This should already be interesting for small ranks, see for example \cite{KhKu-sl3-web-bases,MaPaTu-sl3-web-algebra}, or related quantum base-changes such as \cite{RuTy-specht-web,RuTy-transition-specht-web}. Generating sufficient data for these cases appears feasible, though computational tools for these specific problems are not yet widely available as far as we are aware while writing this paper.

\item \textit{Graded representation theory.}
The simple representations of the symmetric group in finite characteristic have a nontrivial grading by \cite{BrKl-hecke-klr}. By taking \emph{graded dimensions}, one gets a polynomial invariant associated to (certain) partitions of $\sank$. Similar constructions arise in related algebras, such as KLR algebras, {\cf}\cite{KhLa-cat-quantum-sln-first,Ro-2-kac-moody}, graded Temperley--Lieb algebras, {\cf} \cite{PlRyHa-graded-tl-ab}, and KLRW algebras, {\cf} \cite{We-weighted-klr,MaTu-klrw-crystal}. The computational tools in \cite{EvMa-def} may provide a valuable starting point for data generation, with ongoing efforts (as communicated to us by the second author of \cite{EvMa-def}) to integrate these structures into SageMath.

\item \textit{Web spaces invariants.}
Underlying the quantum knot polynomials are spaces of representation theoretical origin, often called \emph{webs} (spiders \cite{Ku-spiders-rank-2} or birdtracks \cite{Cv-bridtracks}), provide a rich setting for studying polynomial invariants. Early examples include \cite{RuTeWe-sl2,Ya-invariant-graphs}. Evaluation of such webs yields polynomials akin to quantum knot polynomials. While the type A exterior case has a closed formula, {\cf} \cite{LaTuVa-web-evaluation}, suggesting limited complexity, extensions to other settings, such as $SO(3)$ invariants and their connections to chromatic polynomials, {\cf} \cite{Ya-invariant-graphs,FeKr-chomatic-category}, remain largely unexplored. The chromatic polynomial has recently been studied via TDA \cite{SaSc-ballmapper}, and the quantum version of this category as, {\eg}, in \cite{MoPeSn-categories-trivalent-vertex,Tu-web-reps} might give some quantization of \cite{SaSc-ballmapper}. Computational complexity here is expected to align with that of canonical bases as in (II). 

\item \textit{Tensor categories and knot/three-manifold invariants.}
Building on \cite{LeHaSa-jones,DlGuSa-mapper}, a promising direction is to analyze knot and three-manifold invariants arising from \emph{tensor categories} (in the sense of \cite{EtGeNiOs-tensor-categories}). Examples include invariants derived from Deligne categories, {\cf} \cite{FlLa-knots}, mixed quantum groups, {\cf} \cite{AnWe-mixed-qgroup,SuTuWeZh-mixed-tilting}, subfactor theory, following the exposition in \cite{Jo-new}, or exotic asymptotic algebras (see, for example, \cite{MaMaMiTuZh-soergel-2reps} for a summary). Currently, there are no well-developed computational approaches in this area as far as we know, making it an exciting frontier for future research.

\item \textit{Many other exciting direction.}
For example, homology and Ext groups naturally yield polynomial invariants, opening another avenue for investigation. Exploring these connections may provide novel insights into representation theory and beyond.

\end{enumerate}
As a side note, many of these and other representation theoretical problems could benefit from AI or deep learning approaches, but it is essential to distinguish these techniques from the big data methodologies emphasized in this work.
\medskip

\noindent\textbf{Acknowledgments.}
We are deeply grateful to Joe Baine, Andrew Mathas, Geordie Williamson, and Victor Zhang for their valuable insights and many stimulating discussions\ochanged{, and like to thank the referee for helpful feedback}.
Special thanks to the organizers of the exceptional conference ``Diagrammatic Intuition and Deep Learning in Mathematics'' (York, July 2024) and to Radmila Sazdanovic, whose inspiring talk at that conference started this project. DT would like to thank the knot PD[X[3, 1, 4, 32], X[1, 7, 2, 6], X[7, 3, 8, 2], X[11, 5, 12, 4], X[5, 11, 6, 10], X[19, 8, 20, 9], X[9, 18, 10, 19], X[25, 12, 26, 13], X[13, 24, 14, 25], X[31, 15, 32, 14], X[15, 21, 16, 20], X[27, 16, 28, 17], X[17, 26, 18, 27], X[21, 29, 22, 28], X[29, 23, 30, 22], X[23, 31, 24, 30]]---you are not just one of many.

AL is grateful for the support and hospitality of the Sydney Mathematical Research Institute (SMRI). AL and PV were supported by a PHC Tournesol Wallonie Bruxelles grant, and DT was sponsored by the ARC Future Fellowship FT230100489. 

In this work, we utilized ChatGPT for coding, proofreading, generating conjectures, and assisting with proof development (and this sentence). Most of the computations were performed using the resources of the Laboratoire de Math{\'e}matiques Blaise Pascal.

Following publication, we were contacted by Eivind Otto Hjelle, who proposed an AI-generated proof of \autoref{C:Growth} (see \autoref{R:ChatGPT} for details). We are very grateful to Eivind for reaching out.
\newpage

\section{\changed{The data and figures}}

\changed{In this section, we have collected all the figures for this paper. The notation is detailed in the main text, but we have crafted the captions and descriptions to be as self-explanatory as possible. Readers are encouraged to browse through them and form their own conjectures.}

\begin{figure}[H]
\fcolorbox{tomato!50}{orchid!10}{
$
\begin{pmatrix}
\fcolorbox{white}{white}{1} & \fcolorbox{white!50}{white!50}{1}
\\
\fcolorbox{black}{black}{{\color{white}0}} & \fcolorbox{white!50}{white!50}{1}
\end{pmatrix}
\leftrightsquigarrow
\begin{tikzpicture}[anchorbase]
\node at (0,0) {\includegraphics[height=3.3cm]{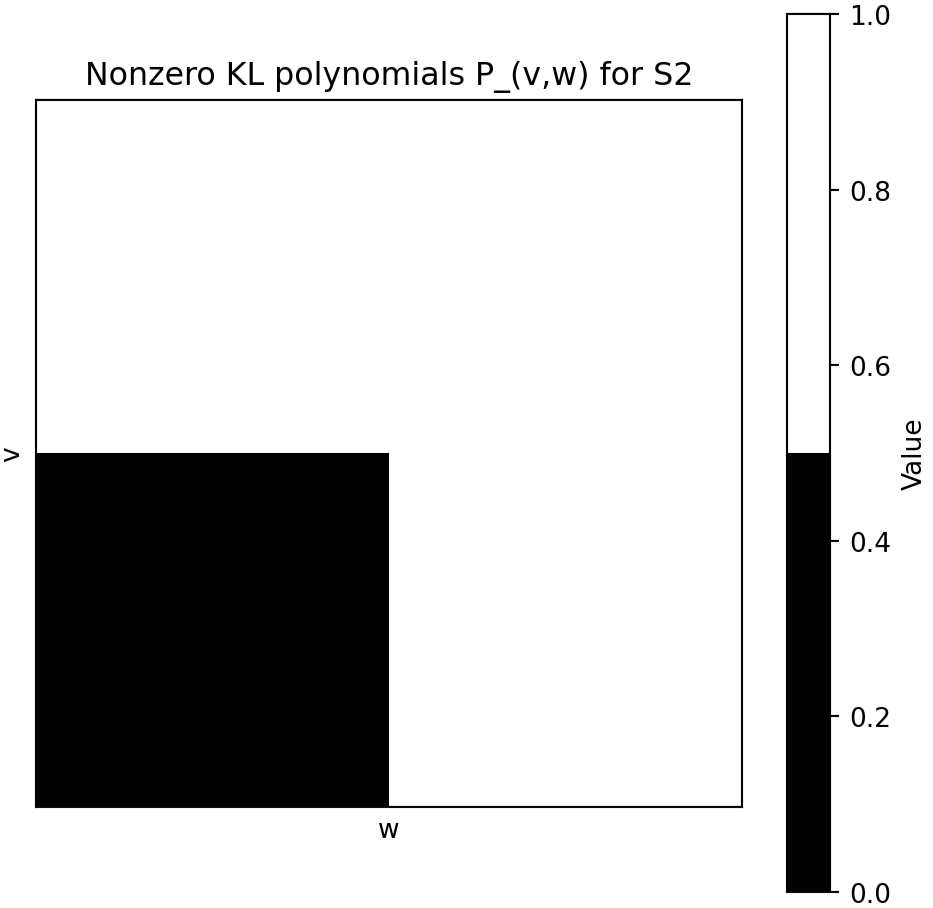}};
\end{tikzpicture}
.
$
}
\caption{Let us explain \emph{grid plots of matrices}: One visualizes a matrix with $\R_{\geq 0}$ entries as is with a heat map, where the temperature corresponds to the size of the entry. We use black as the coldest=smallest entry and white as the hottest=biggest entry. $H$ is the matrix with heat plot and $B$ its binary version. This figure shows an example of a matrix $B=H$.}
\label{sec3-first}
\end{figure}

\begin{figure}[H]
\fcolorbox{tomato!50}{orchid!10}{\includegraphics[height=5.2cm]{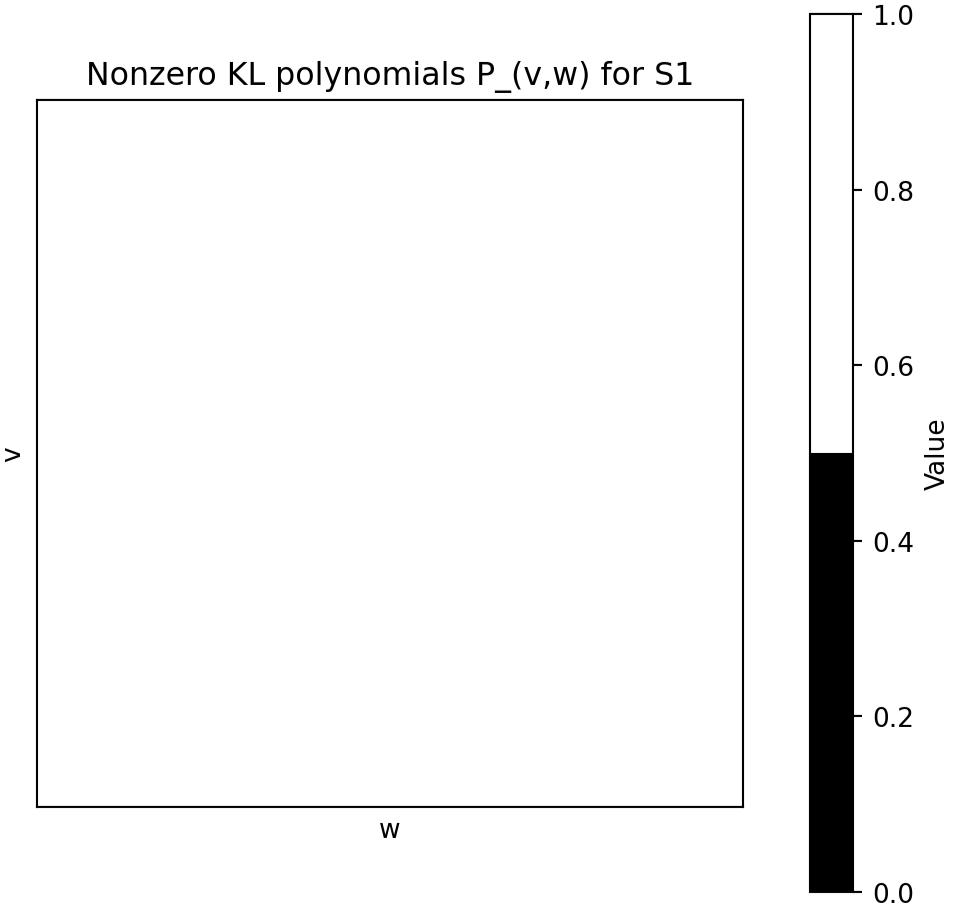}\includegraphics[height=5.2cm]{figs/kla1.png}\includegraphics[height=5.2cm]{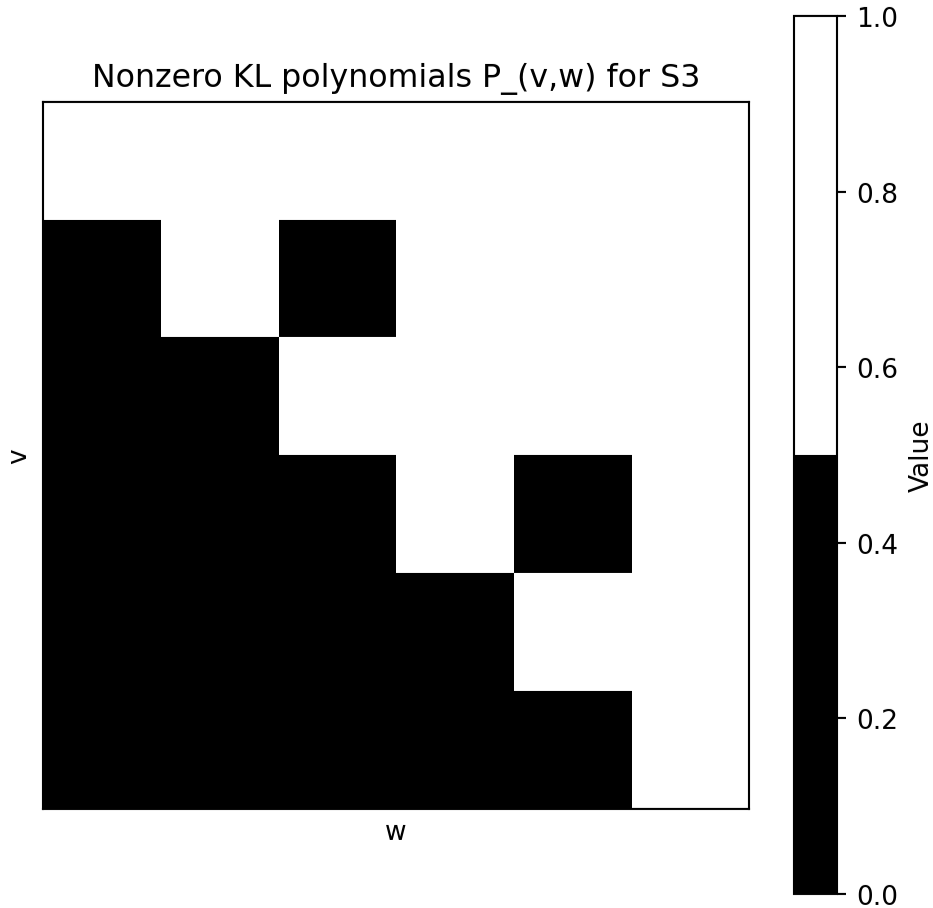}}
\caption{The plot of the matrices $B=H$ for evaluations of KL polynomials at $\vpar=1$, with conventions as in \autoref{sec3-first}, for $\sank\in\{1,2,3\}$.}
\label{sec3-1}
\end{figure}

\begin{figure}[H]
\fcolorbox{tomato!50}{orchid!10}{\includegraphics[height=6.0cm]{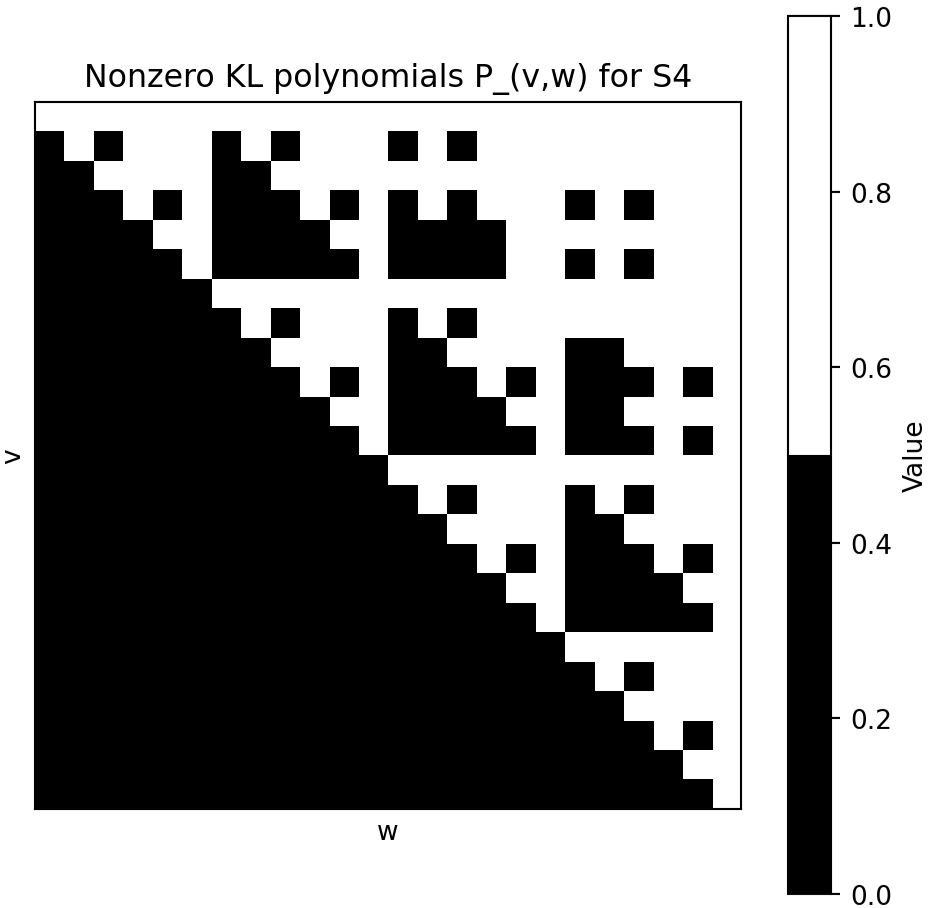}\includegraphics[height=6.0cm]{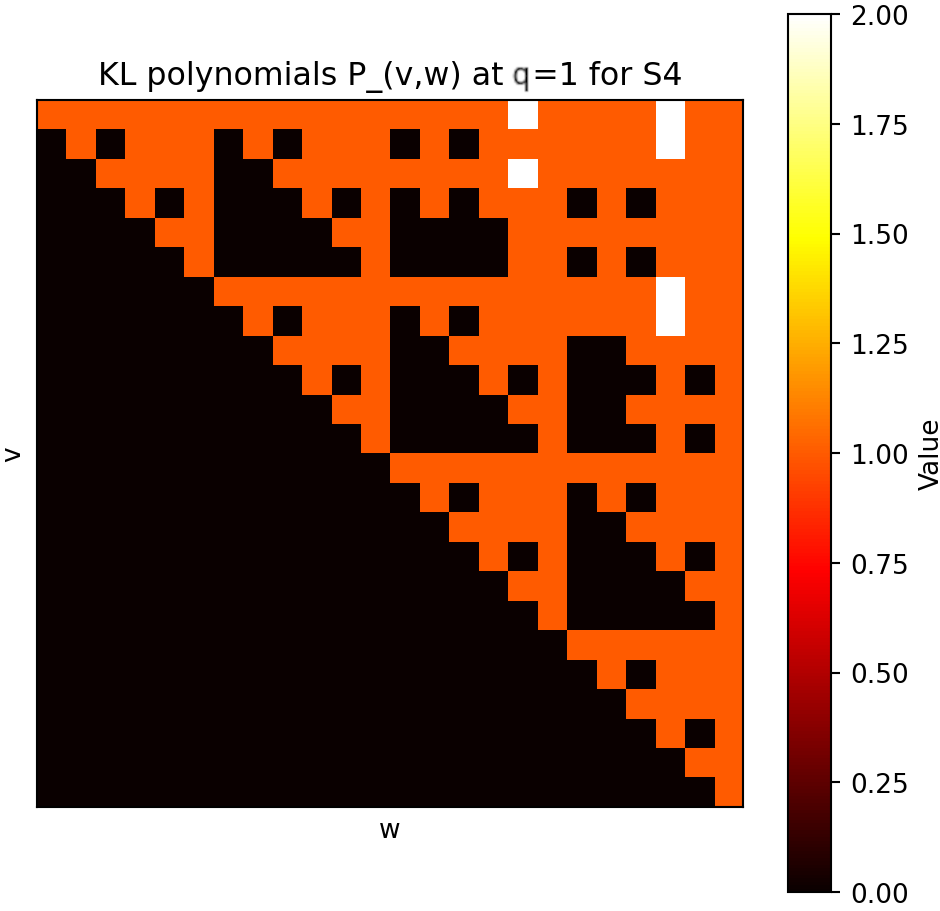}}
\fcolorbox{tomato!50}{orchid!10}{\includegraphics[height=6.0cm]{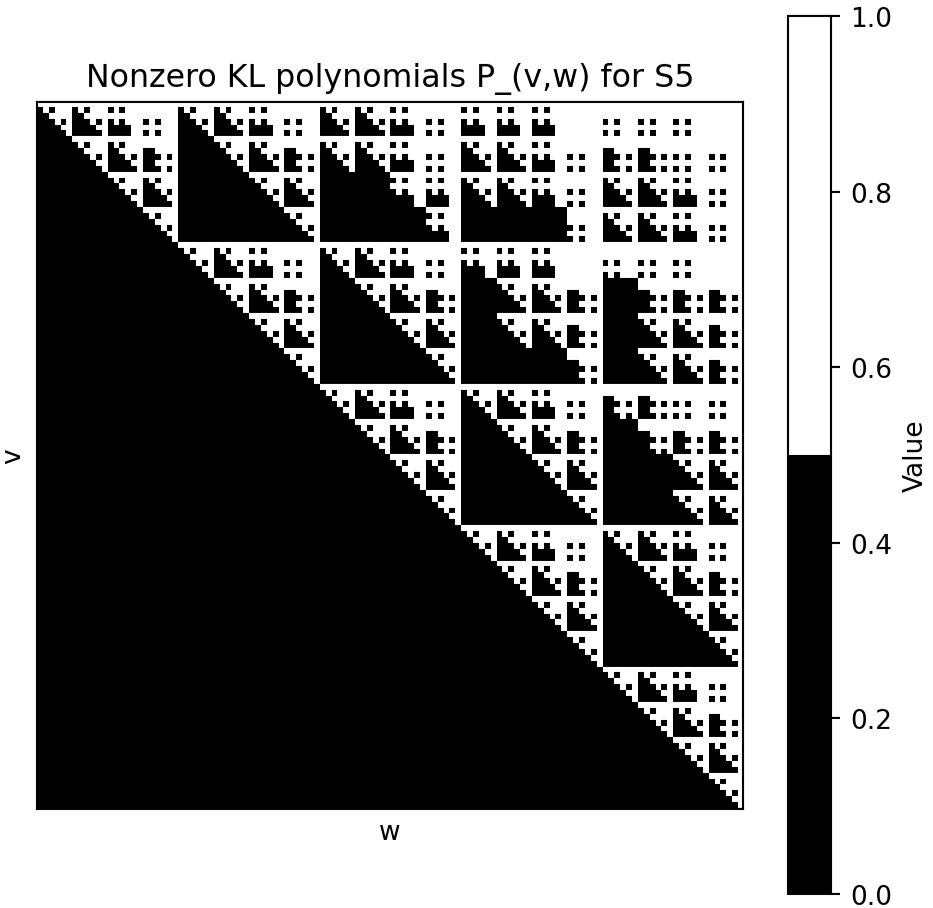}\includegraphics[height=6.0cm]{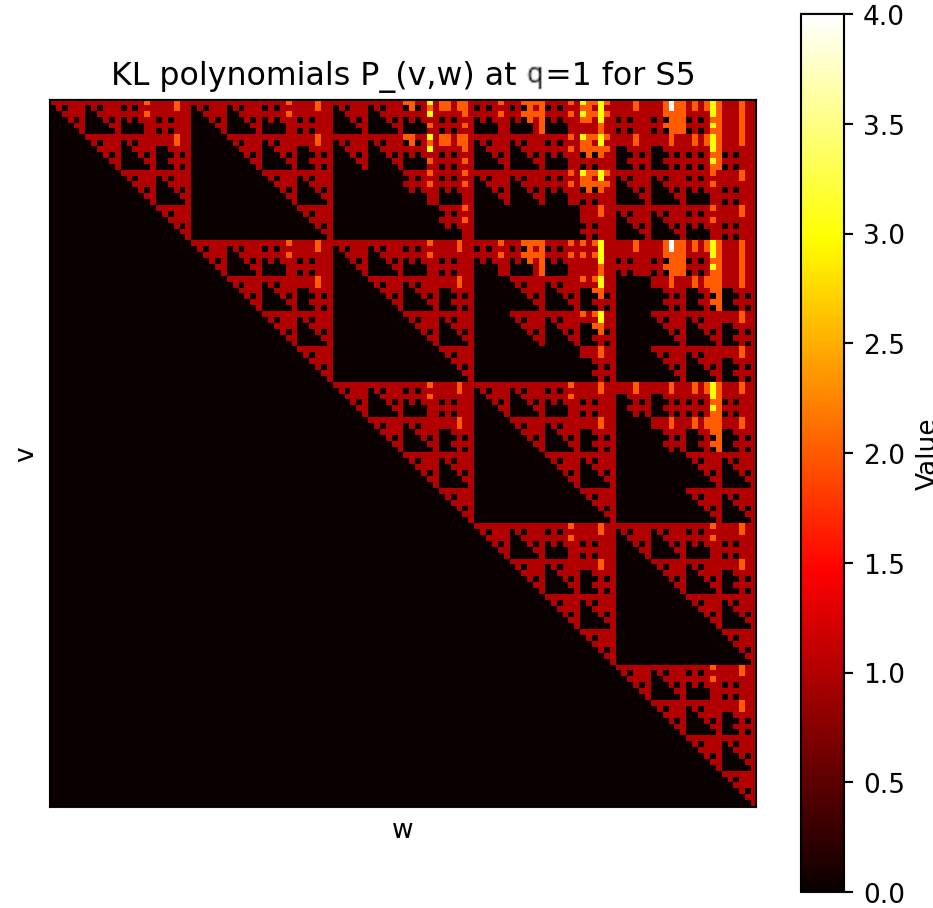}}
\fcolorbox{tomato!50}{orchid!10}{\includegraphics[height=6.0cm]{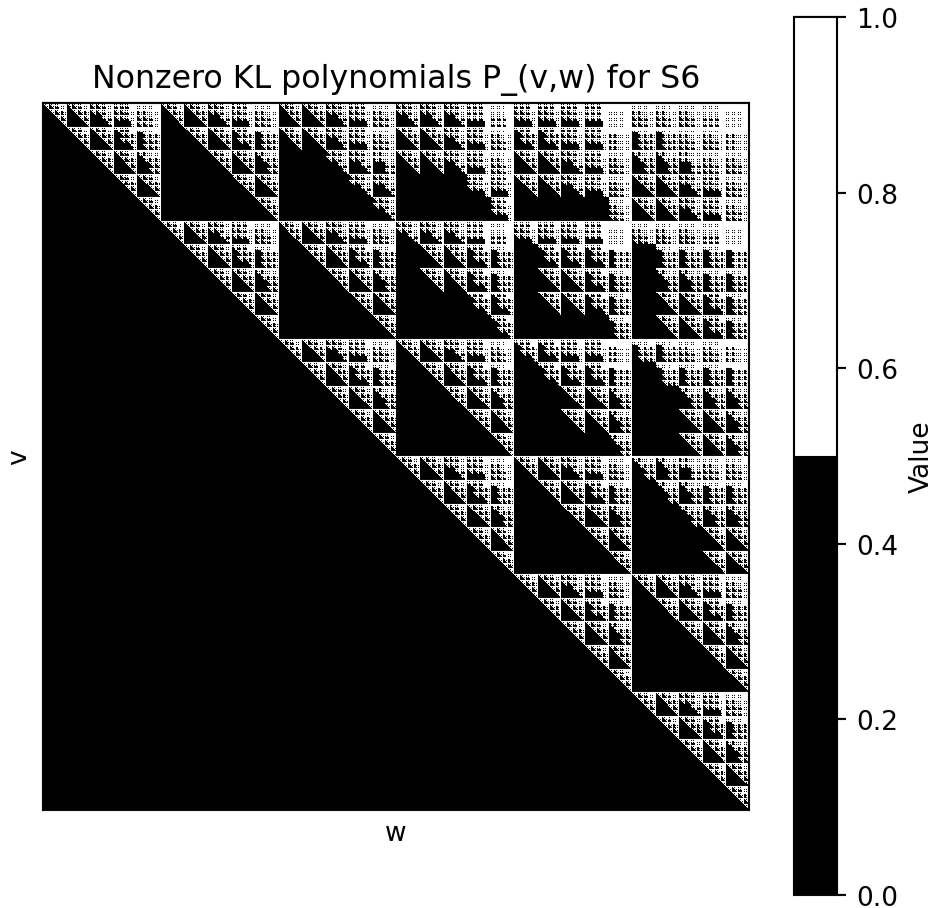}\includegraphics[height=6.0cm]{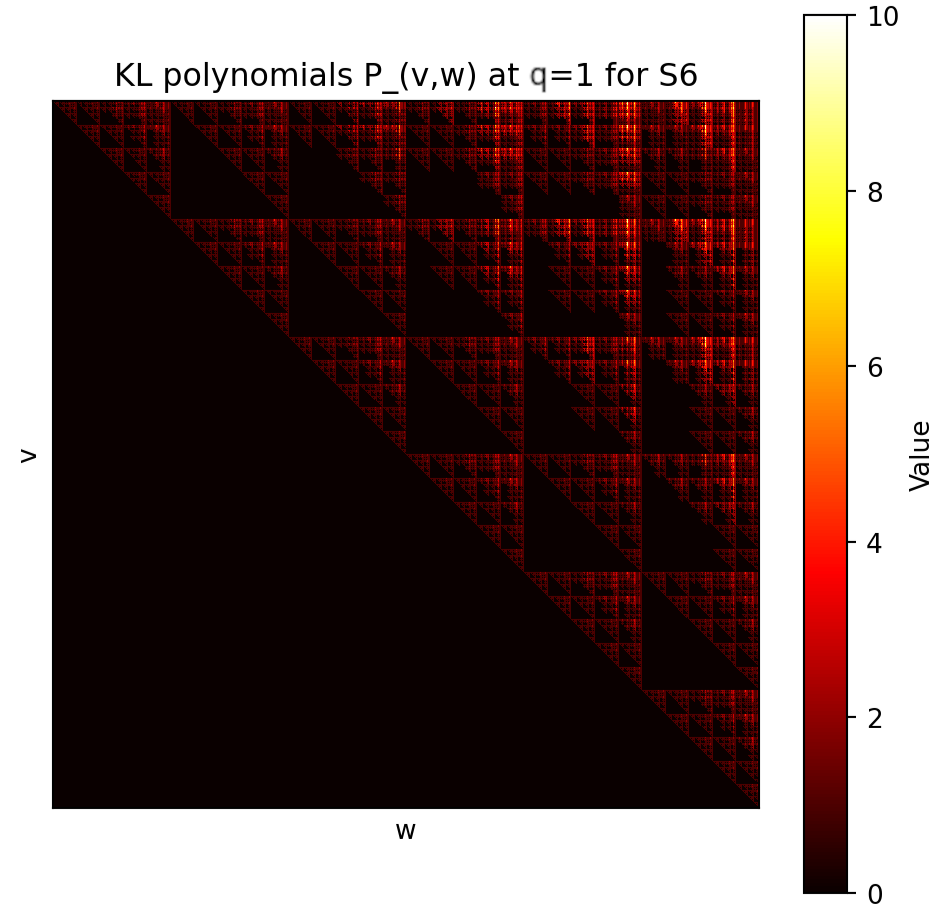}}
\fcolorbox{tomato!50}{orchid!10}{\includegraphics[height=6.0cm]{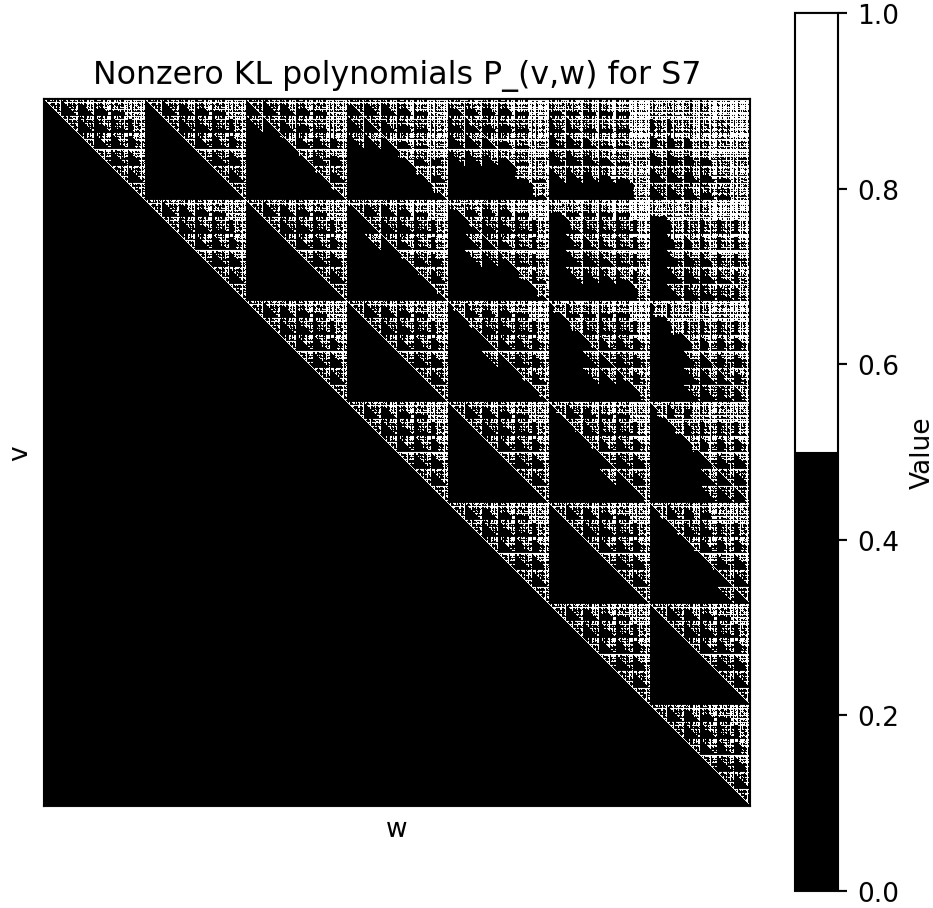}\includegraphics[height=6.0cm]{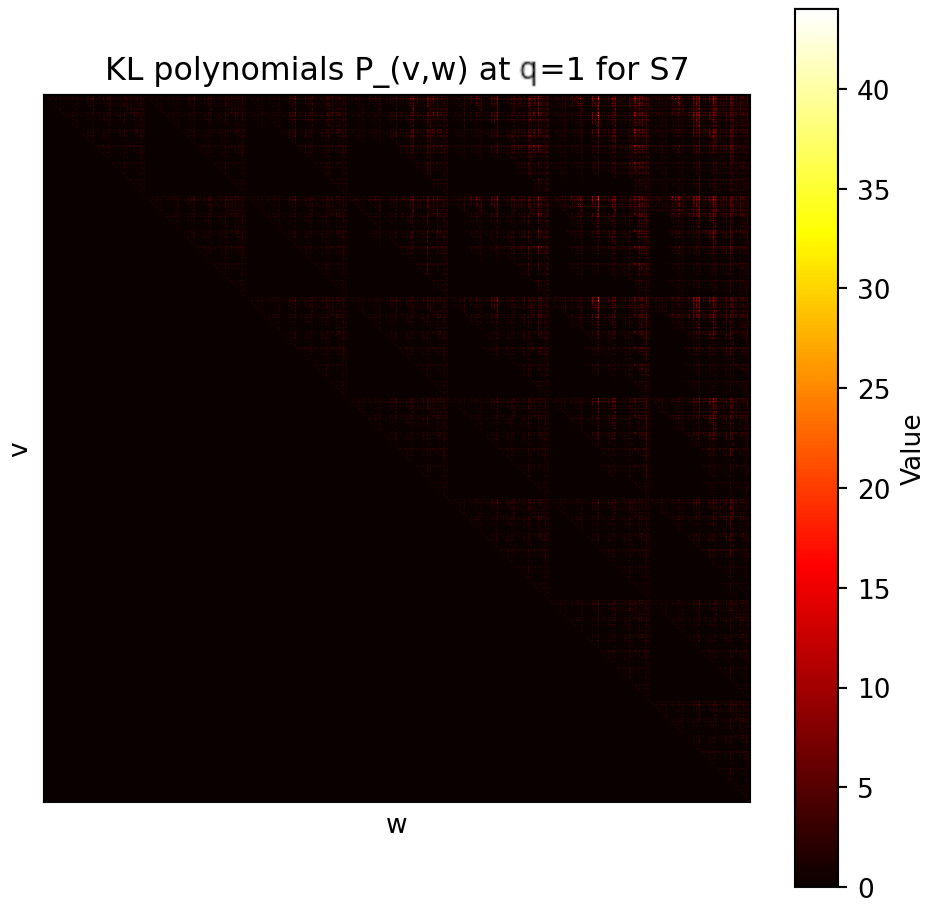}}
\caption{The plot as in \autoref{sec3-1} but for $\sank\in\{4,5,6,7\}$.}
\label{sec3-2}
\end{figure}

\begin{figure}[H]
\fcolorbox{tomato!50}{orchid!10}{\includegraphics[height=4.2cm]{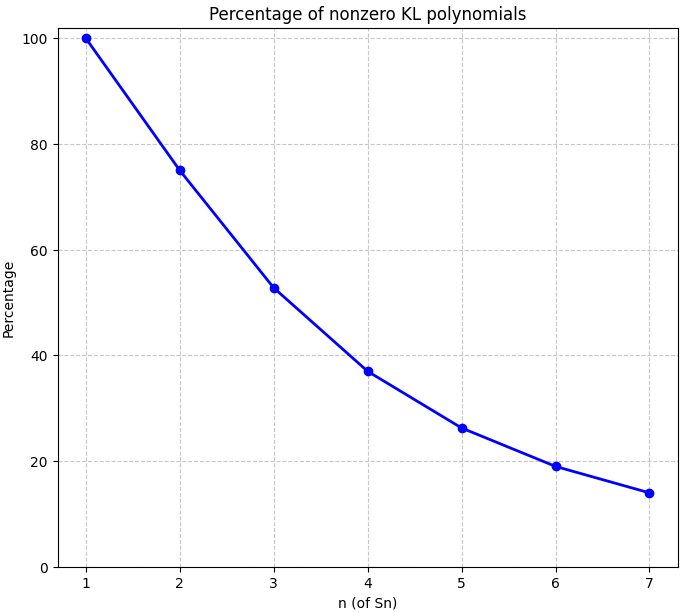}\includegraphics[height=4.2cm]{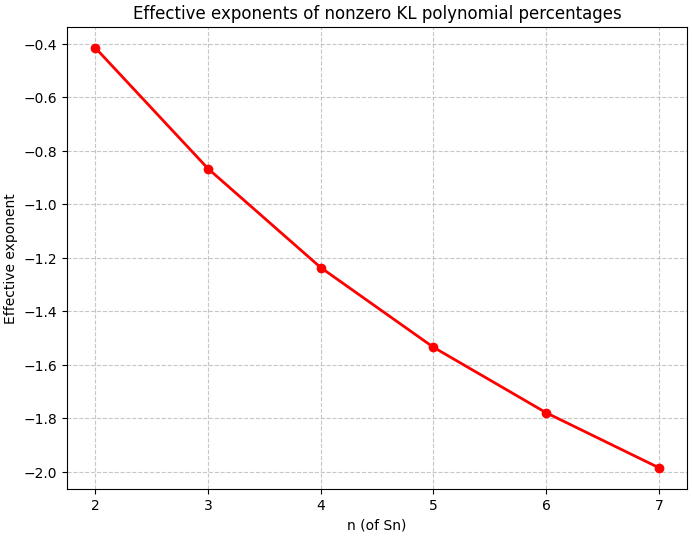}}
\fcolorbox{tomato!50}{orchid!10}{\includegraphics[height=4.2cm]{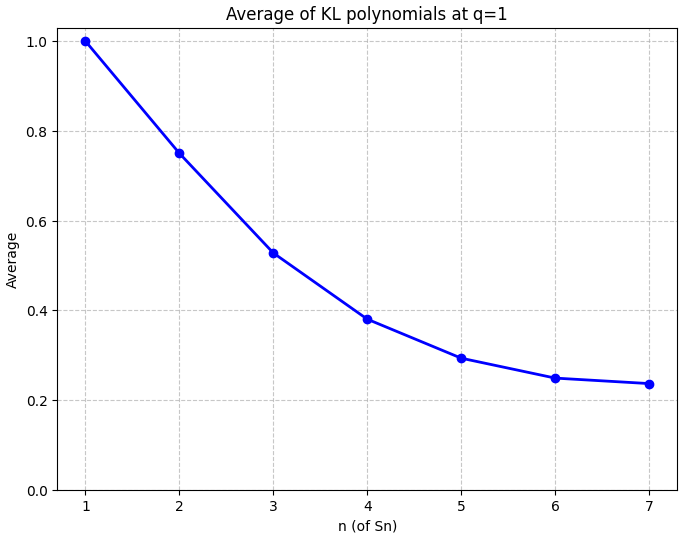}\includegraphics[height=4.2cm]{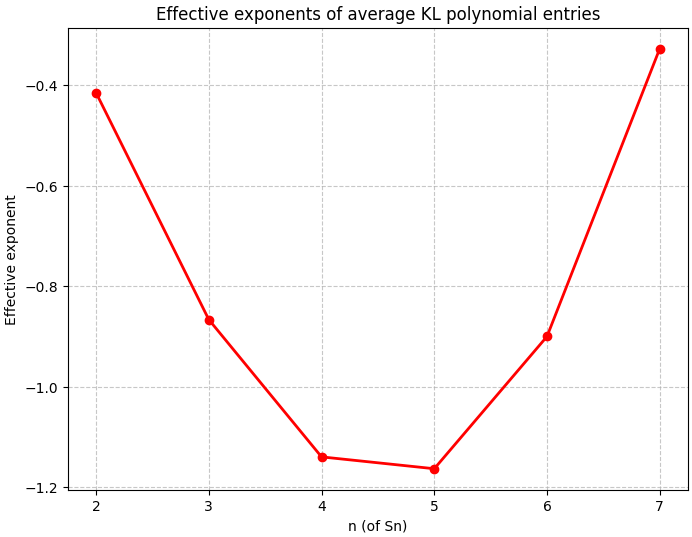}}
\fcolorbox{tomato!50}{orchid!10}{\includegraphics[height=4.2cm]{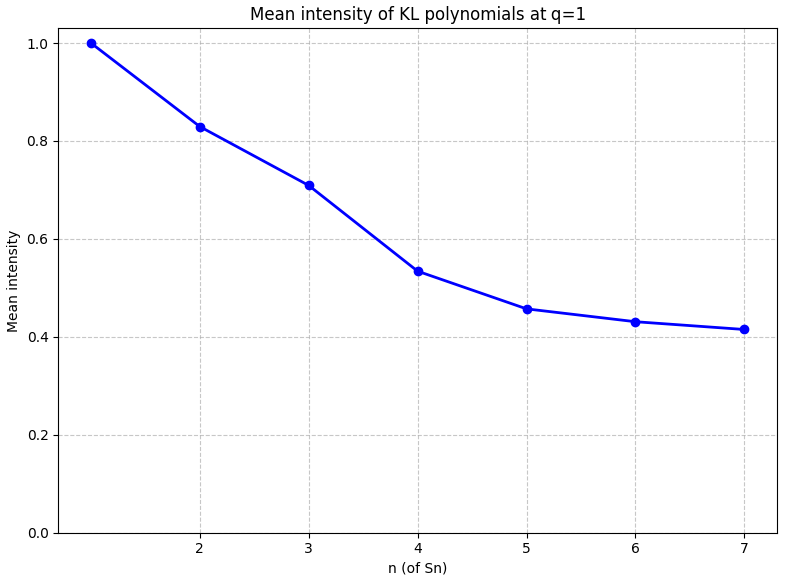}\includegraphics[height=4.2cm]{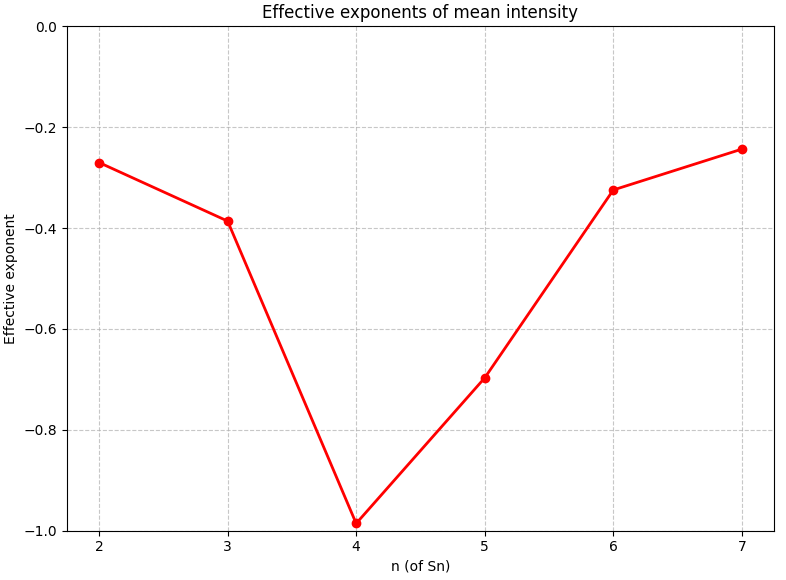}}
\begin{gather*}
\begin{tabular}{c|ccccccccccc}
$\sank$  & 1 & 2 & 3 & 4 & 5 & 6 & 7 \\
\hline
\% $\neq 0$ & 100 & 75 & 52.78 & 36.98 & 26.26 & 18.98 & 13.98 \\
\hline
av & 1 & 0.75 & 0.5278 & 0.3802 & 0.2933 & 0.2489 & 0.2366 \\
\hline
in & 1 & 0.8295 & 0.7095 & 0.5346 & 0.4573 & 0.4307 & 0.4140 \\
\end{tabular}
\end{gather*}
\caption{Plots of various statistics associated to the percentage of nonzero KL polynomials. From top to bottom: the percentage of nonzero polynomials, the average evaluation and the intensity of the heat plots. The right plots illustrate the approximate $a\in\R$ assuming that the function on the left grows like $\sank^{a}$. Note that the average value and the percentage of nonzero  entries are the same for $\sank\in\{1,2,3\}$.}
\label{sec3-3}
\end{figure}
\vspace{-0.5cm}
\begin{figure}[H]
\fcolorbox{tomato!50}{orchid!10}{\begin{tikzpicture}[anchorbase]
\node at (0,0) {\includegraphics[height=5.0cm]{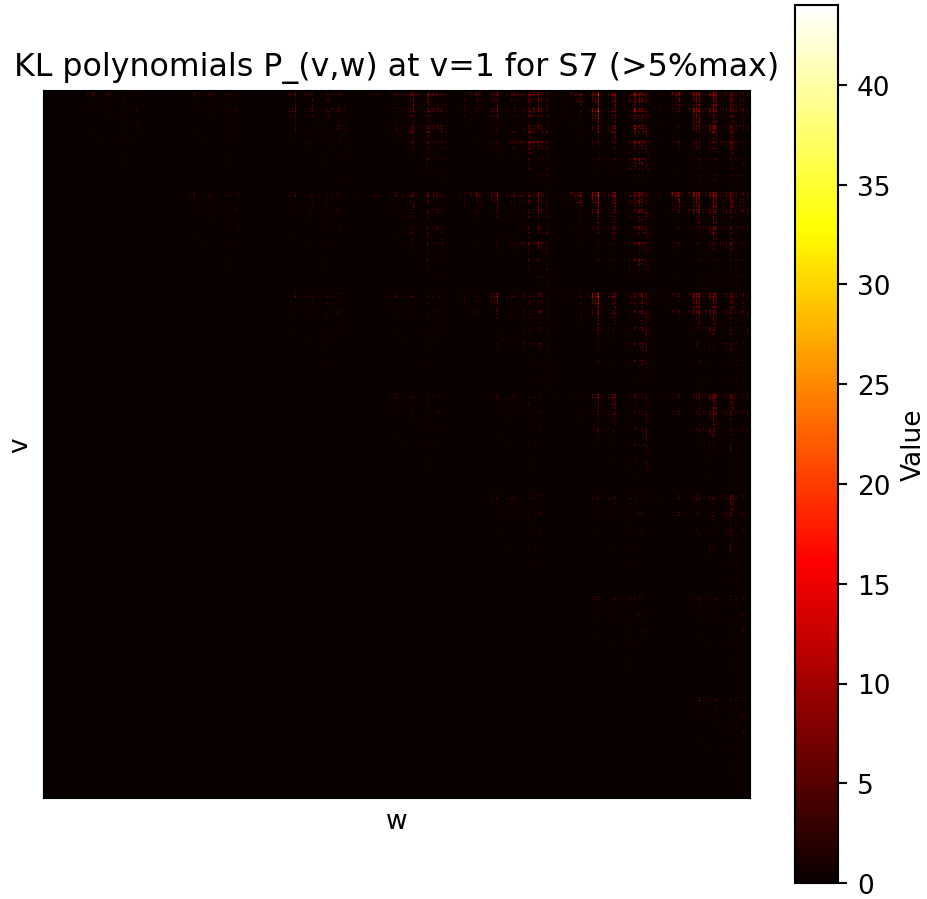}};
\end{tikzpicture}
,\hspace{-0.2cm}
\begin{tikzpicture}[anchorbase]
\node at (0,0) {\includegraphics[height=5.0cm]{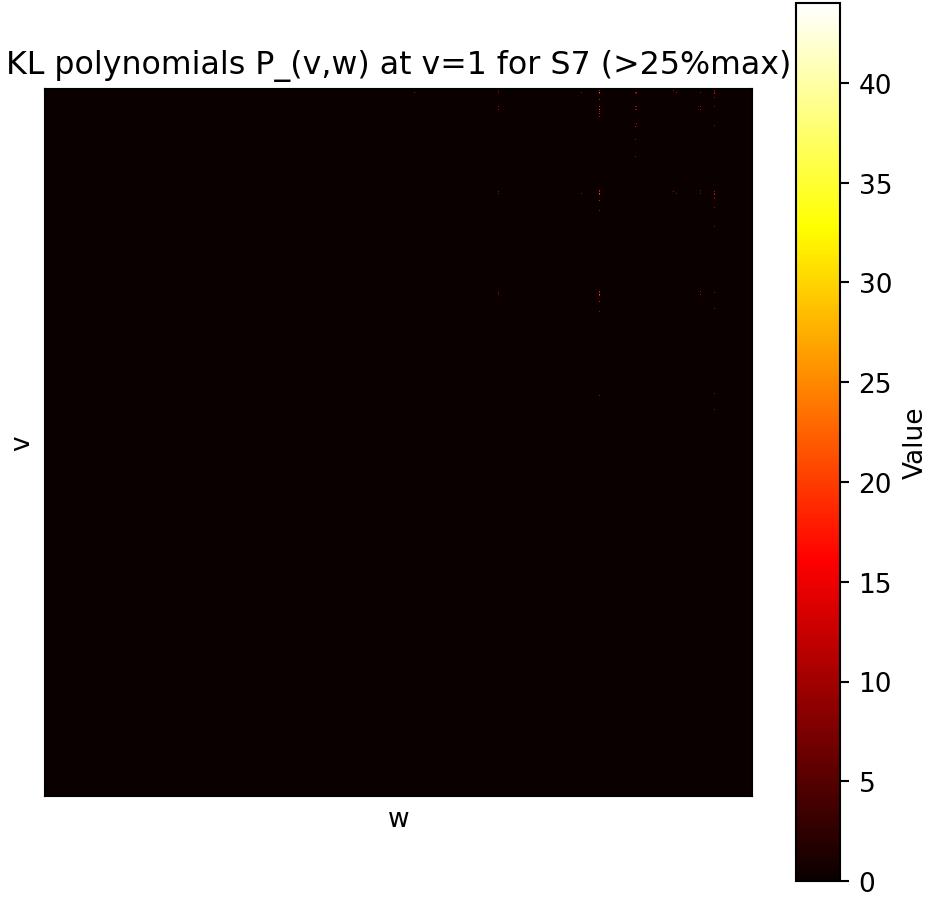}};
\end{tikzpicture}
,\hspace{-0.2cm}
\begin{tikzpicture}[anchorbase]
\node at (0,0) {\includegraphics[height=5.0cm]{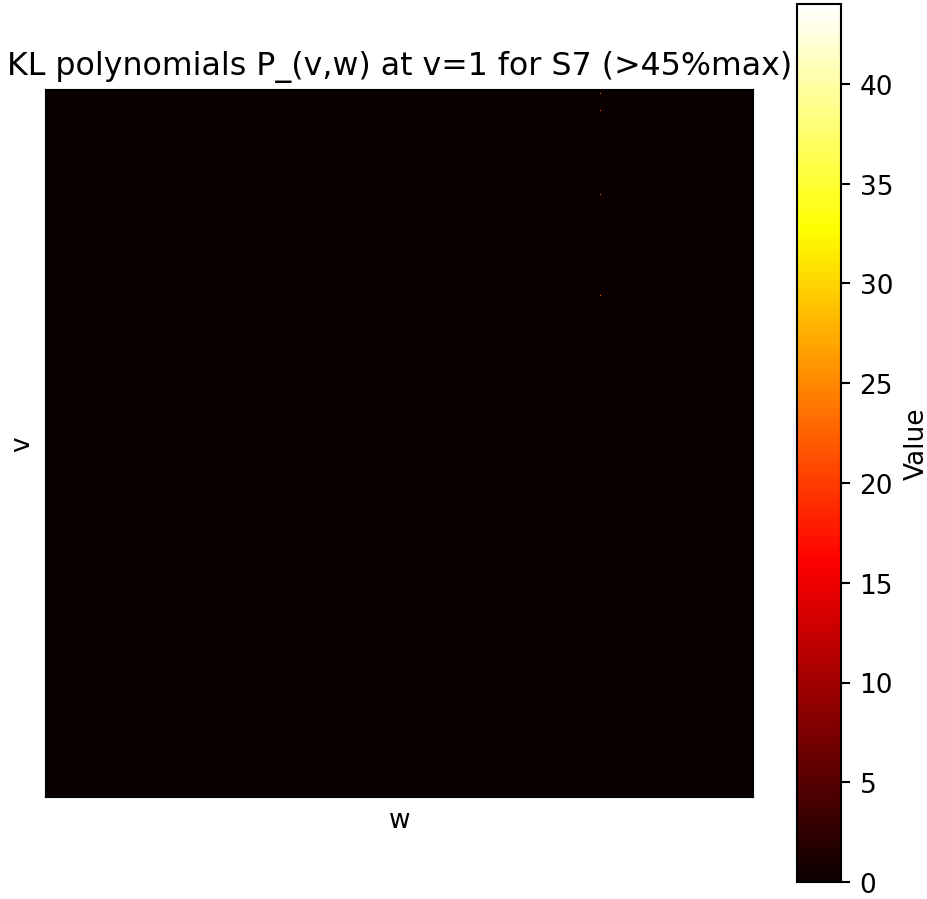}};
\end{tikzpicture}
}
\caption{The pictures, from left to right, show the pairs $v,w\in\sym[7]$ such that 
$\KL{v}{w}(1)$ is greater than $x$\% of $\max\{\KL{v}{w}(1)|v,w\in\sym[7]\}$ for $x=5,25,45$, respectively. The rare red and white pixels can better be seen on the high resolution plots that can be found on \cite{LaTuVa-kl-big-data-code}.}
\label{sec3-4}
\end{figure}

\begin{figure}[H]
\fcolorbox{tomato!50}{orchid!10}{\includegraphics[height=5.5cm]{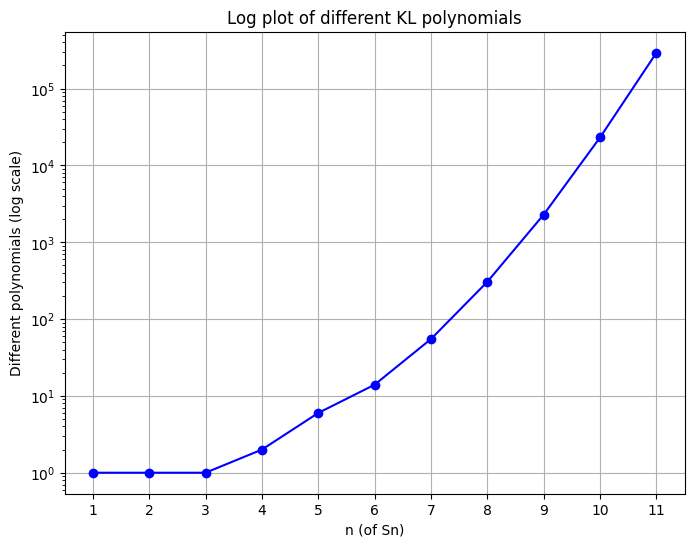}\includegraphics[height=5.5cm]{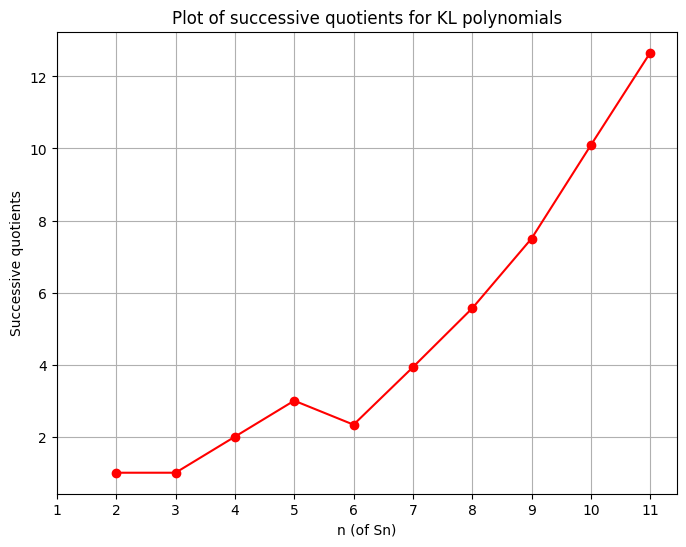}}
\begin{gather*}
\begin{tabular}{c|ccccccccccc}
$\sank$  & 1 & 2 & 3 & 4 & 5 & 6 & 7 & 8 & 9 & 10 & 11 \\
\hline
\# & 1 & 1 & 1 & 2 & 6 & 14 & 55 & 306 & 2295 & 23163 & 293189 \\
\hline
\% & 100 & 50 & 16.67 & 8.33 & 5 & 1.94 & 1.09 & 0.76 & 0.63 & 0.64 & 0.73 \\
\end{tabular}
\end{gather*}
\caption{The number of different KL polynomials $\kl$ and their percentage. The graph on the left is a log plot. The right plot illustrates the approximate $a\in\R$ assuming that the function on the left grows like $a^{\sank}$.}
\label{sec4}
\end{figure}

\begin{figure}[H]
\fcolorbox{tomato!50}{orchid!10}{\includegraphics[height=5.5cm]{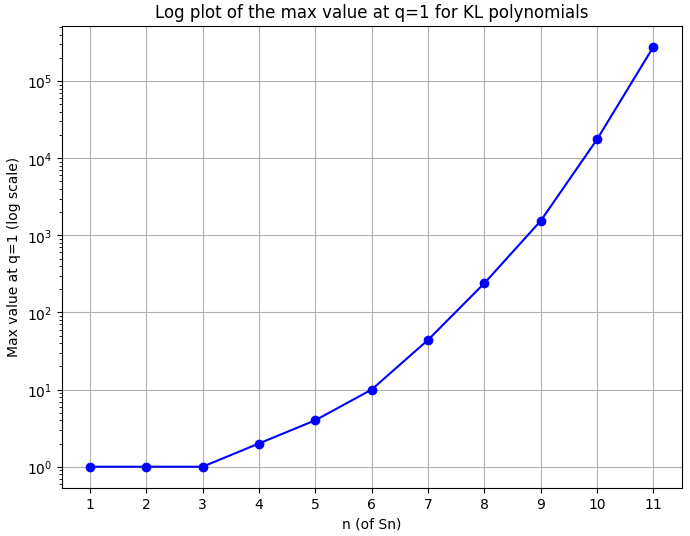}\includegraphics[height=5.5cm]{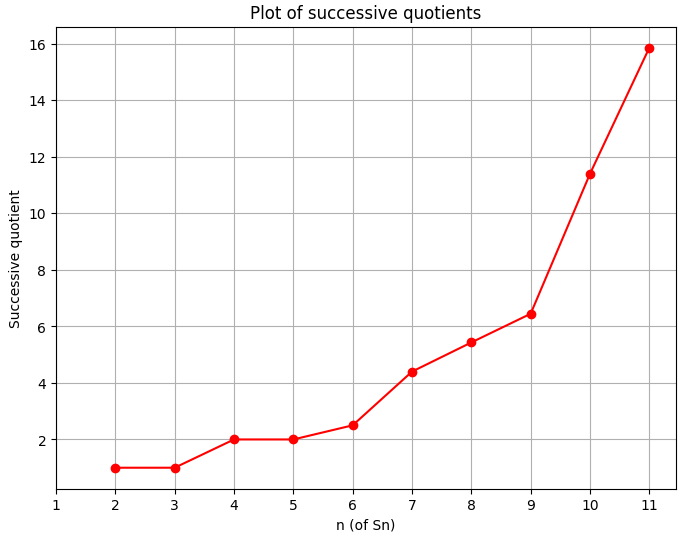}}
\fcolorbox{tomato!50}{orchid!10}{\includegraphics[height=5.5cm]{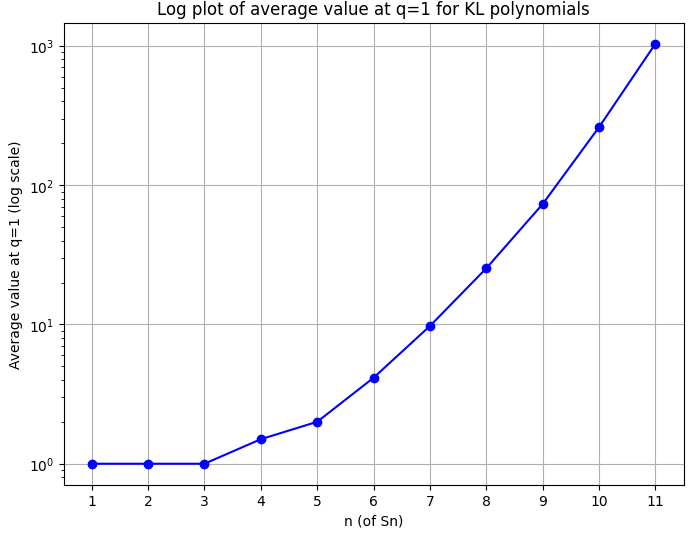}\includegraphics[height=5.5cm]{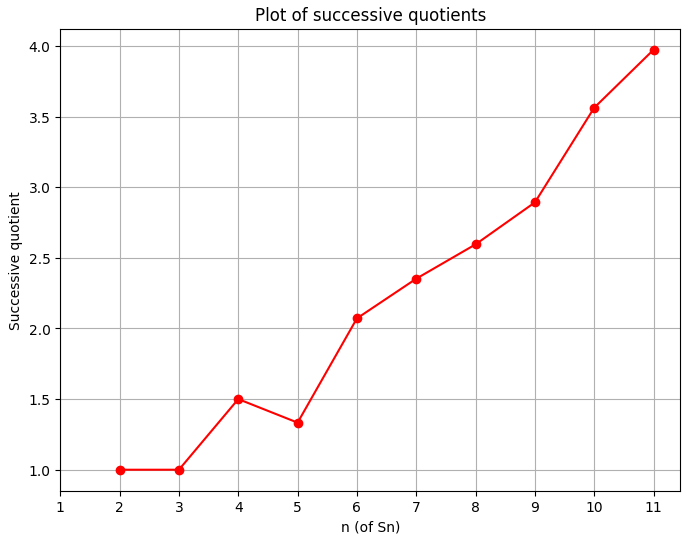}}
\begin{gather*}
\begin{tabular}{c|ccccccccccc}
$\sank$  & 1 & 2 & 3 & 4 & 5 & 6 & 7 & 8 & 9 & 10 & 11 \\
\hline
max & 1 & 1 & 1 & 2 & 4 & 10 & 44 & 239 & 1541 & 17566 & 278576 \\
\hline
av & 1 & 1 & 1 & 1.5 & 2.4 & 4.85 & 10 & 25.30 & 72.71 & 259.03 & 1035.78 \\
\end{tabular}
.
\end{gather*}
\caption{The maximal and average value of the evaluation at $\vpar=1$. The graphs on the left are log plots. The right plots illustrate the approximate $a\in\R$ assuming that the function on the left grows like $a^{\sank}$.}
\label{sec5-evaluation}
\end{figure}

\begin{figure}[H]
\fcolorbox{tomato!50}{orchid!10}{\includegraphics[height=5.5cm]{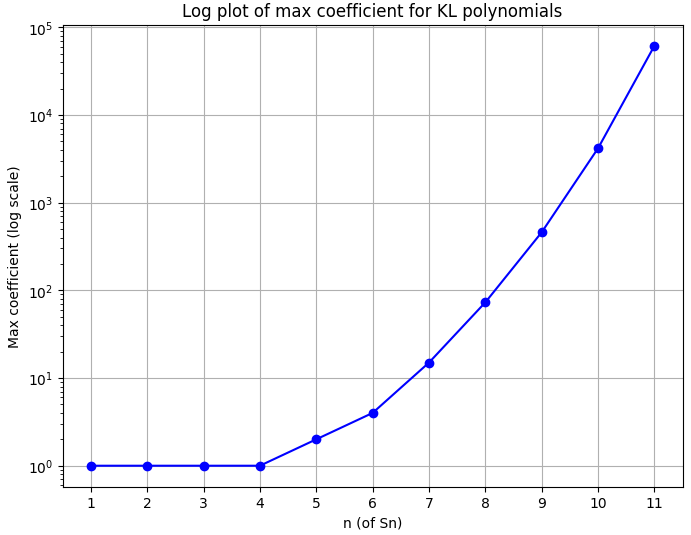}\includegraphics[height=5.5cm]{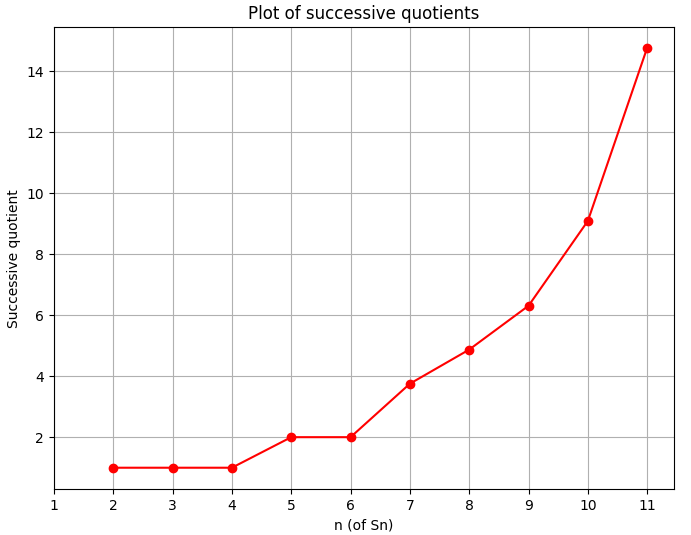}}
\fcolorbox{tomato!50}{orchid!10}{\includegraphics[height=5.5cm]{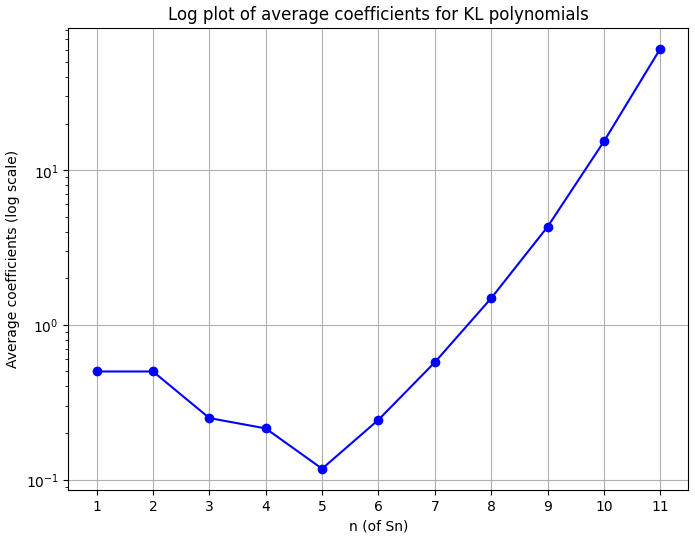}\includegraphics[height=5.5cm]{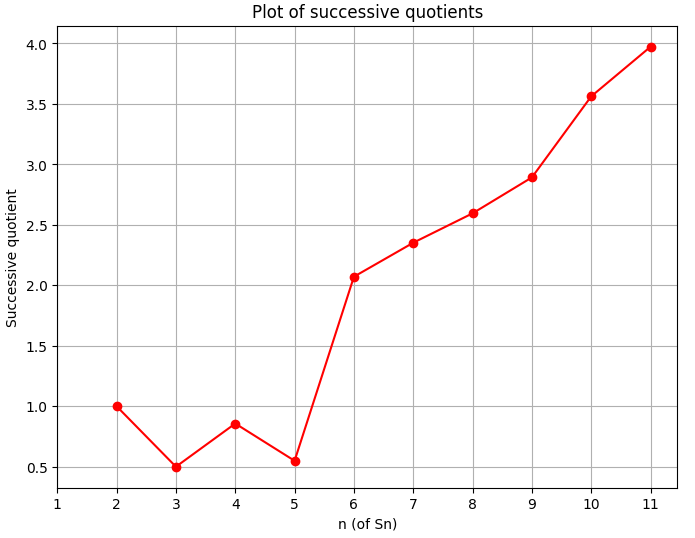}}
\begin{gather*}
\begin{tabular}{c|ccccccccccc}
$\sank$  & 1 & 2 & 3 & 4 & 5 & 6 & 7 & 8 & 9 & 10 & 11 \\
\hline
max & 1 & 1 & 1 & 1 & 2 & 4 & 15 & 73 & 460 & 4176 & 61582 \\
\hline
av & 1 & 0.5 & 0.25 & 0.21 & 0.12 & 0.24 & 0.57 & 1.49 & 4.30 & 15.33 & 39.84 \\
\end{tabular}
.
\end{gather*}
\caption{The maximal and average value of the maximal coefficient. The graphs on the left are log plots. The right plots illustrate the approximate $a\in\R$ assuming that the function on the left grows like $a^{\sank}$.}
\label{sec5-maximum}
\end{figure}

\begin{figure}[H]
\fcolorbox{tomato!50}{orchid!10}{\includegraphics[height=5.5cm]{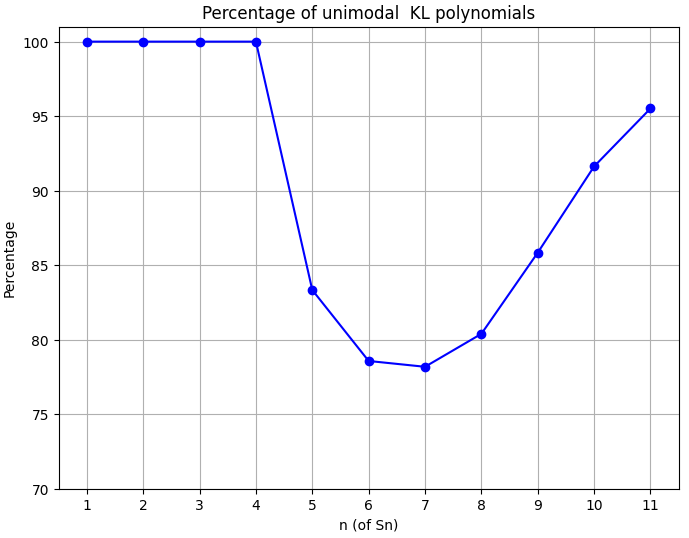}\includegraphics[height=5.5cm]{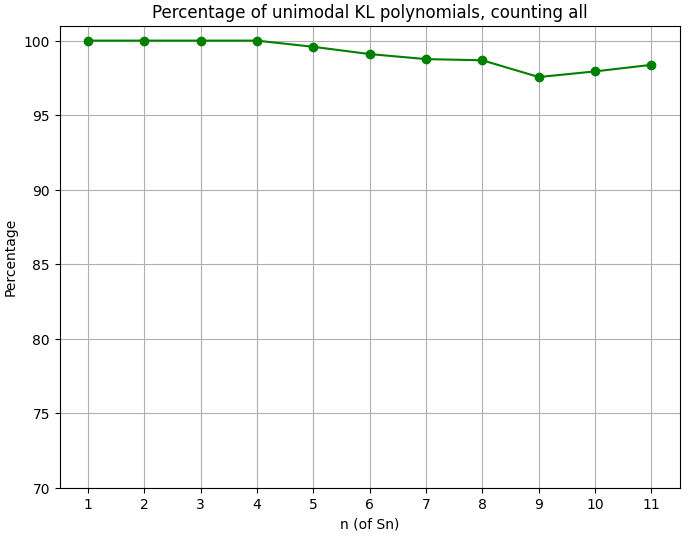}}
\begin{gather*}
\begin{tabular}{c|ccccccccccc}
$\sank$  & 1 & 2 & 3 & 4 & 5 & 6 & 7 & 8 & 9 & 10 & 11 \\
\hline
\% uni(set) & 100 & 100 & 100 & 100 & 83.33 & 78.57 & 78.18 & 80.39 & 85.84 & 91.62 & 95.51 \\
\hline
\% uni(multiset) & 100 & 100 & 100 & 100 & 99.59 & 99.10 & 98.76 & 98.68 & 97.56 & 97.93 & 98.38 \\
\end{tabular}
.
\end{gather*}
\caption{Percentage of unimodal (one peak) KL polynomials, counted with multiplicity on the left and without on the right. Note that the bottom value is 70\% not 0\%.}
\label{sec6-unimodal}
\end{figure}

\begin{figure}[H]
\fcolorbox{tomato!50}{orchid!10}{\includegraphics[height=5.4cm]{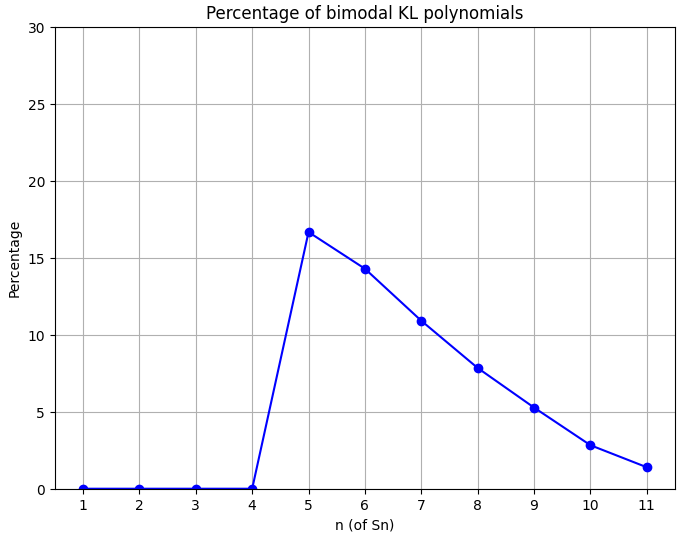}\includegraphics[height=5.4cm]{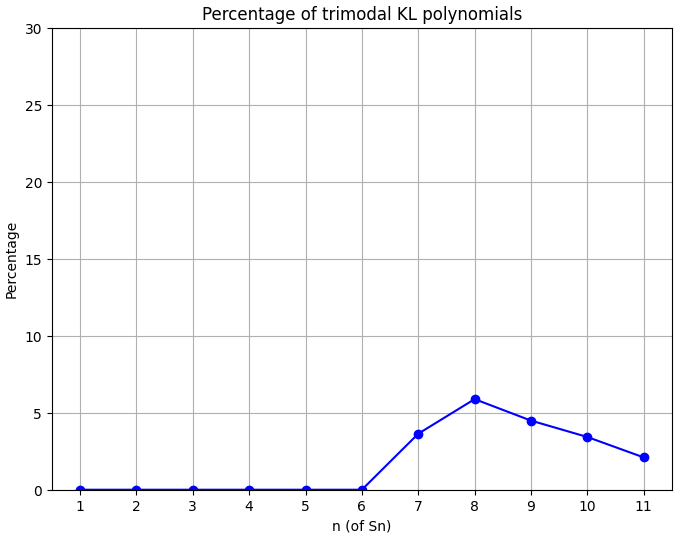}}
\begin{gather*}
\begin{tabular}{c|ccccccccccc}
$\sank$  & 1 & 2 & 3 & 4 & 5 & 6 & 7 & 8 & 9 & 10 & 11 \\
\hline
\% bi & 0 & 0 & 0 & 0 & 16.67 & 14.29 & 10.91 & 7.84 & 5.27 & 2.83 & 1.40 \\
\hline
\% tri & 0 & 0 & 0 & 0 & 0 & 0 & 3.64 & 5.88 & 4.49 & 3.42 & 2.10 \\
\end{tabular}
\end{gather*}
\caption{Percentage of bimodal (two peaks) and trimodal (three peaks) KL polynomials. Note that the top value is 30\% not 100\%.}
\label{sec6-multimodal}
\end{figure}

\begin{figure}[H]
\fcolorbox{tomato!50}{orchid!10}{\includegraphics[height=5.4cm]{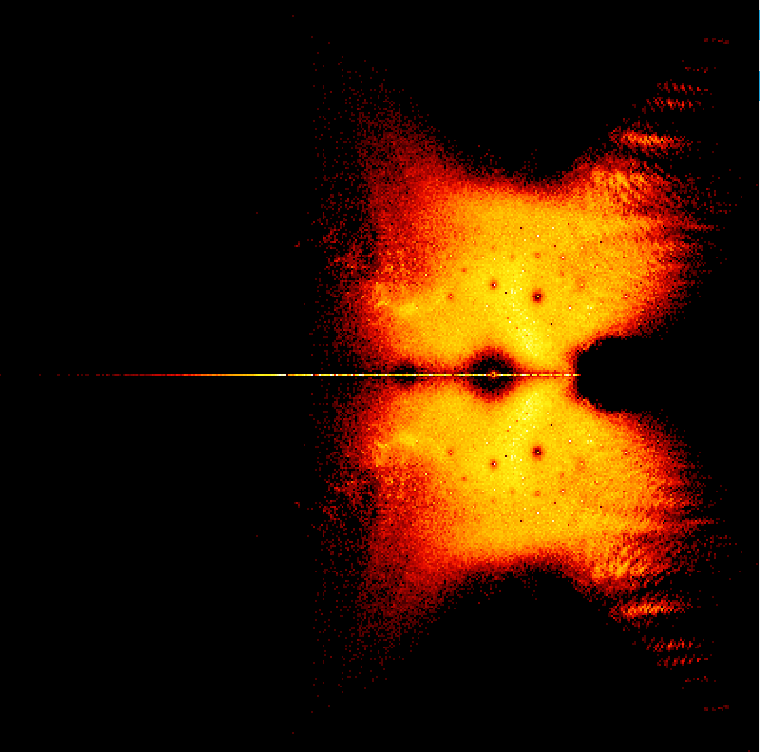}\includegraphics[height=5.4cm]{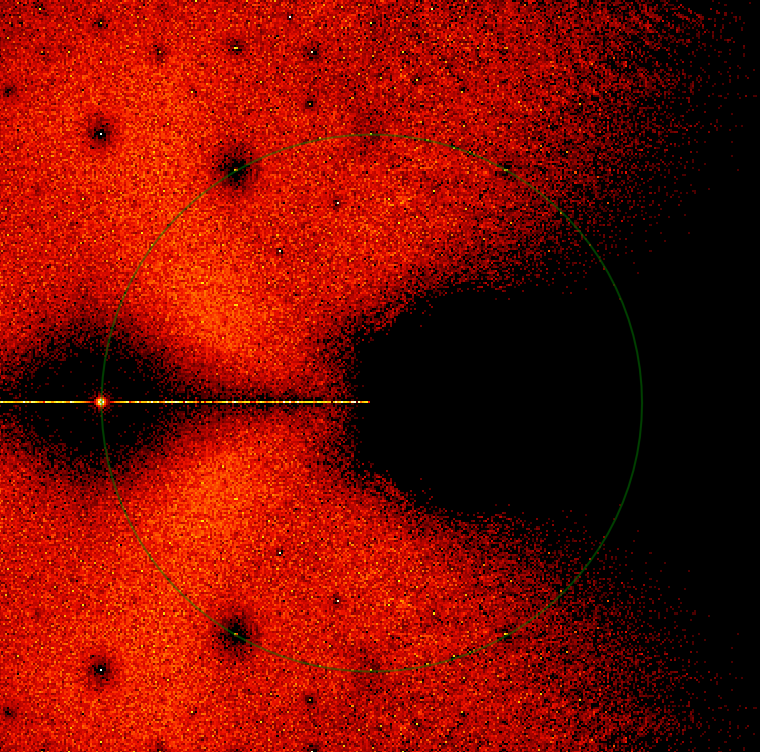}}
\caption{Roots of KL polynomials (counted with multiplicity and illustrated in a heat plot) with the leftmost points being $\approx -6.6532$ and $-1.5$, {\ie} the right picture is a zoom into the whole picture.}
\label{sec7-KL}
\end{figure}

\begin{figure}[H]
\fcolorbox{tomato!50}{orchid!10}{\includegraphics[height=5.4cm]{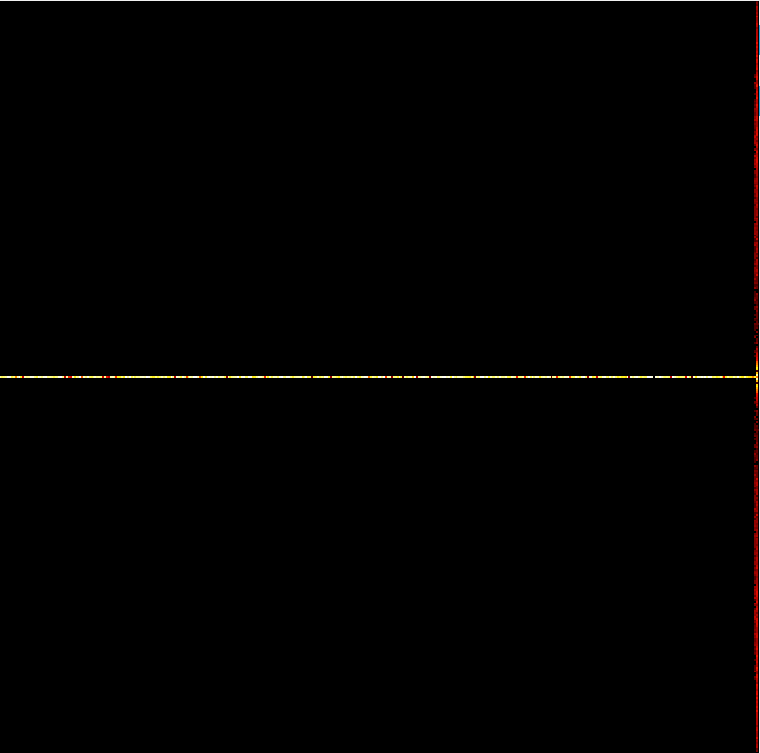}\includegraphics[height=5.4cm]{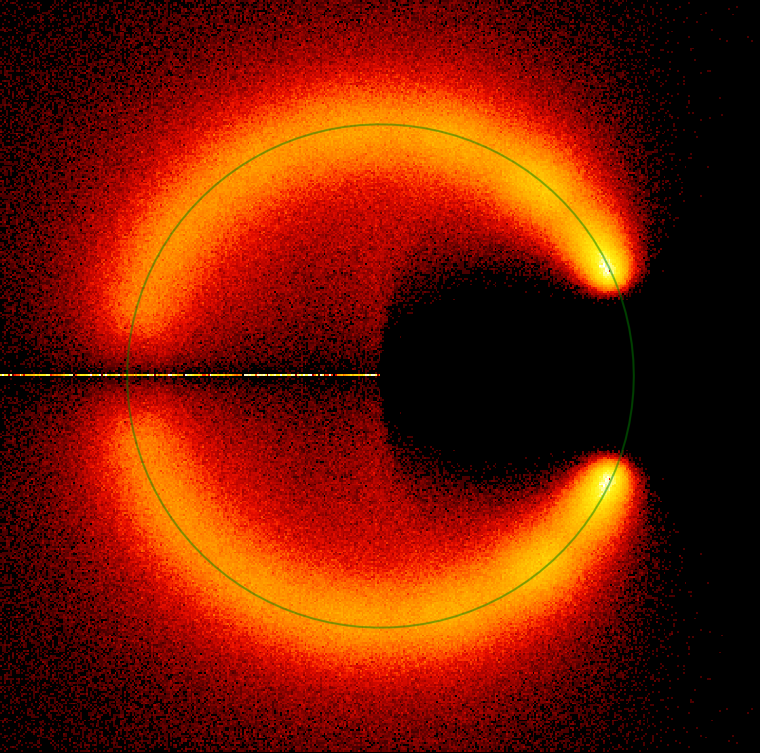}}
\caption{Roots of random polynomials (counted with multiplicity and illustrated in a heat plot) with leftmost points being $\approx -61581.9988$ and $-1.5$, {\ie} the right picture is a zoom into the whole picture. The percentage of real roots is $\approx 24.12\%$, the percentage of roots with|root| in $[0.9, 1.1]$ is $\approx 44.15\%$, the average real value is $\approx -3849.98$.}
\label{sec7-random}
\end{figure}

\begin{figure}[H]
\fcolorbox{tomato!50}{orchid!10}{\includegraphics[height=5.2cm]{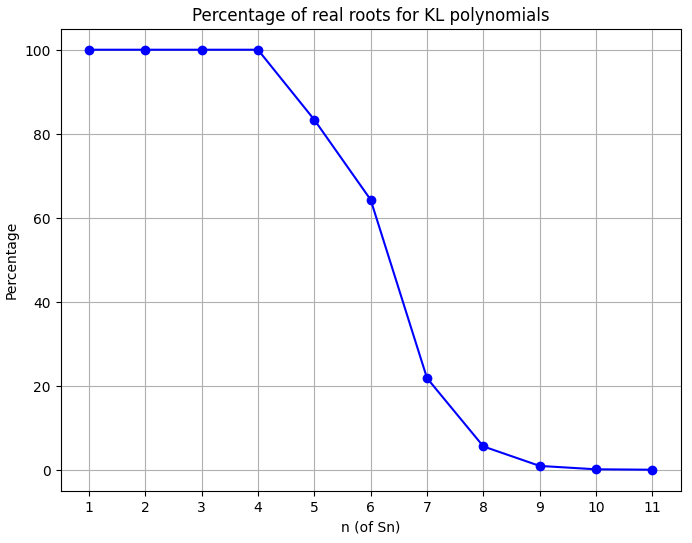}\includegraphics[height=5.2cm]{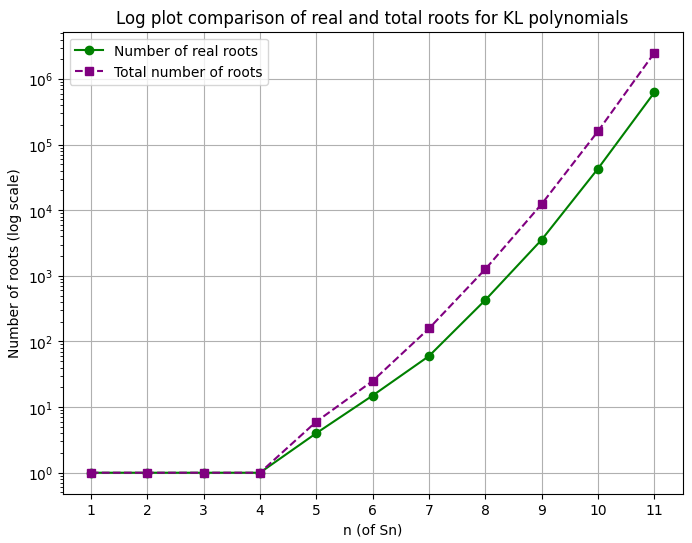}}
\begin{gather*}
\begin{tabular}{c|ccccccccccc}
$\sank$  & 1 & 2 & 3 & 4 & 5 & 6 & 7 & 8 & 9 & 10 & 11 \\
\hline
\% & 100 & 100 & 100 & 100 & 83.33 & 64.29 & 21.82 & 5.56 & 0.92 & 0.11 & 0.01 \\
\end{tabular}
.
\end{gather*}
\caption{Percentage of real roots of KL polynomials.}
\label{sec7-real}
\end{figure}
\vspace{-0.5cm}
\begin{figure}[H]
\fcolorbox{tomato!50}{orchid!10}{\includegraphics[height=5.2cm]{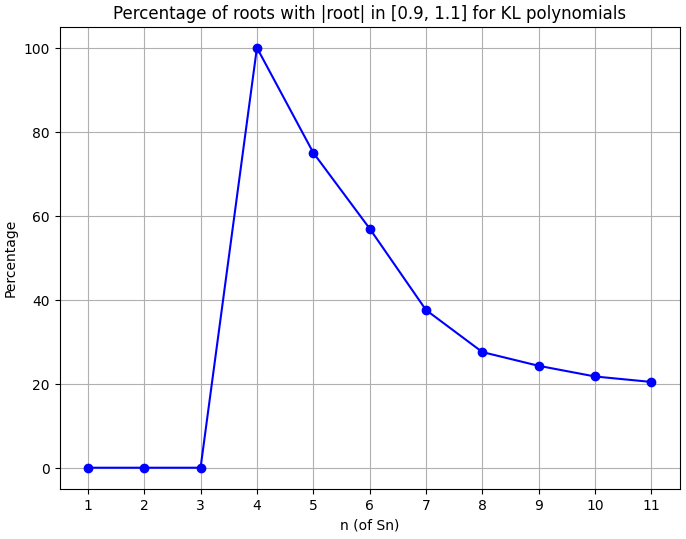}\includegraphics[height=5.2cm]{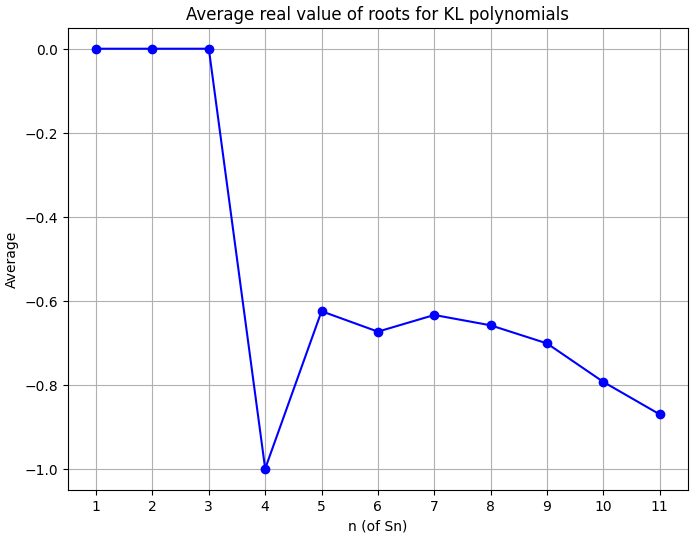}}
\begin{gather*}
\begin{tabular}{c|ccccccccccc}
$\sank$  & 1 & 2 & 3 & 4 & 5 & 6 & 7 & 8 & 9 & 10 & 11 \\
\hline
\% & 0 & 0 & 0 & 100 & 75 & 56.94 & 37.55 & 27.53 & 24.26 & 21.72 & 20.42 \\
av. real & 0 & 0 & 0 & -1 & -0.625 & -0.6736 & -0.6339 & -0.6585 & -0.7016 & -0.7933 & -0.8706 \\
\end{tabular}
.
\end{gather*}
\caption{Percentage of roots around the unit circle and average real value of roots of KL polynomials.}
\label{sec7-distance}
\end{figure}

\begin{figure}[H]
\fcolorbox{tomato!50}{orchid!10}{\includegraphics[height=5.2cm]{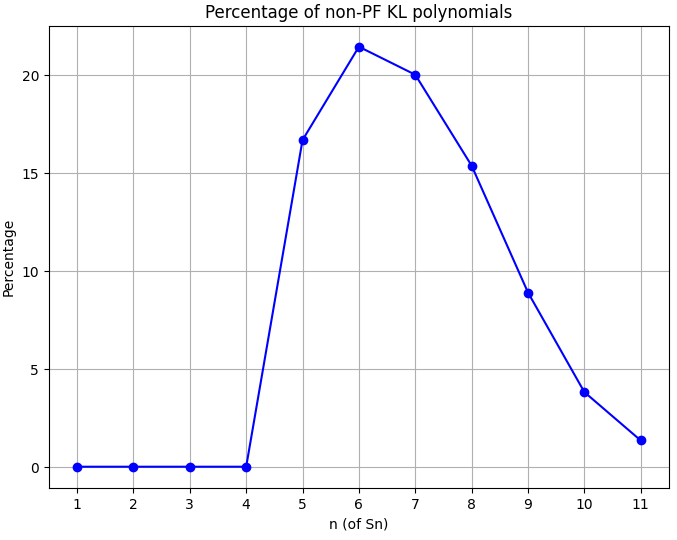}\includegraphics[height=5.2cm]{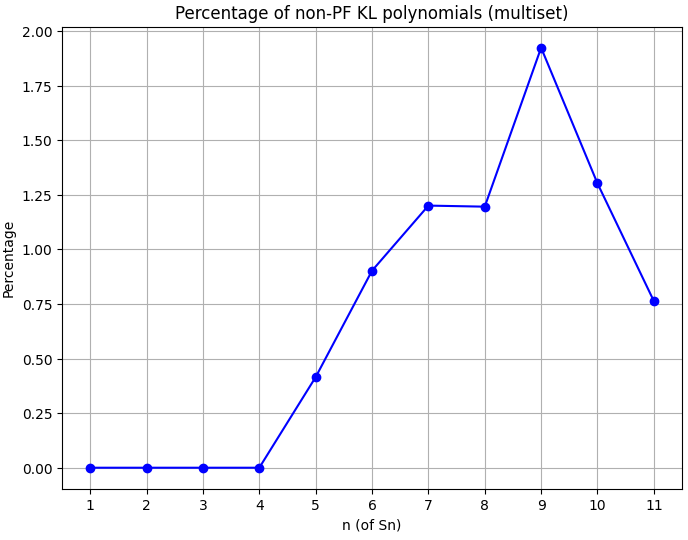}}
\begin{gather*}
\begin{tabular}{c|ccccccccccc}
$\sank$  & 1 & 2 & 3 & 4 & 5 & 6 & 7 & 8 & 9 & 10 & 11 \\
\hline
\% set & 0 & 0 & 0 & 0 & 16.67 & 21.43 & 20 & 15.36 & 8.89 & 3.81 & 1.34 \\
\hline
\% multiset & 0 & 0 & 0 & 0 & 0.41 & 0.90 & 1.2 & 1.20 & 1.92 & 1.30 & 0.76 \\
\end{tabular}
.
\end{gather*}
\caption{Percentage of KL polynomials not satisfying the PF property (they do not have a leading real root), counted with multiplicity on the left and without on the right. Note that the top of diagrams are not at 100\%.}
\label{sec8-PF}
\end{figure}

\begin{figure}
\fcolorbox{tomato!50}{orchid!10}{\includegraphics[height=5.2cm]{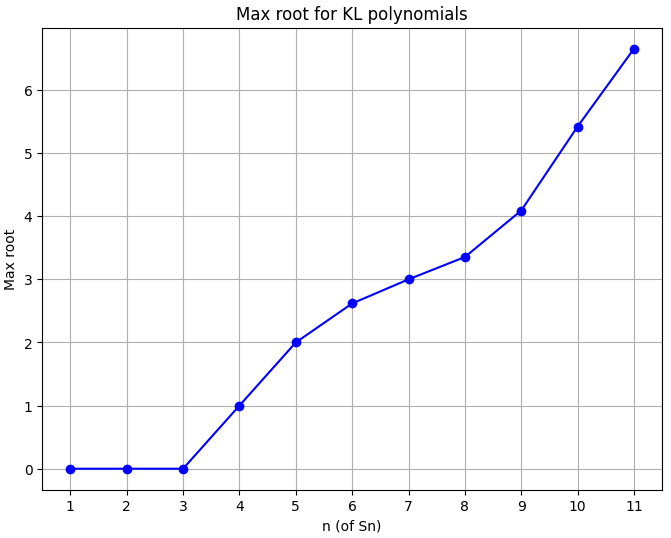}\includegraphics[height=5.2cm]{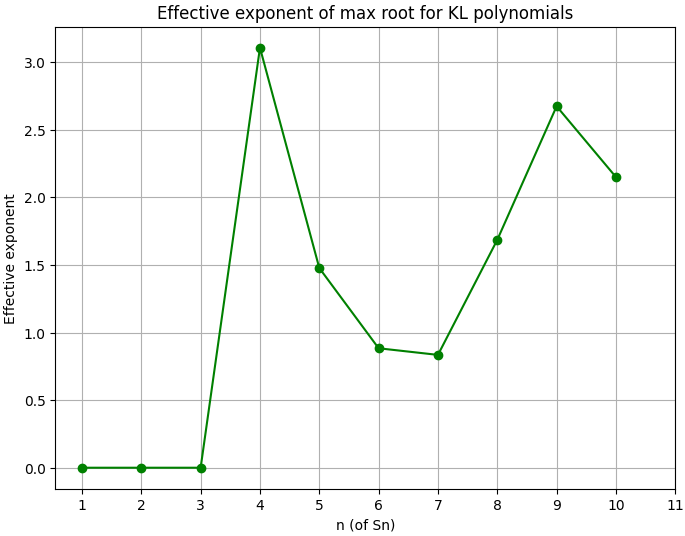}}
\begin{gather*}
\begin{tabular}{c|ccccccccccc}
$\sank$  & 1 & 2 & 3 & 4 & 5 & 6 & 7 & 8 & 9 & 10 & 11 \\
\hline
max & 0 & 0 & 0 & 1 & 2 & 2.618 & 3 & 3.354 & 4.090 & 5.4206 &  6.6532 \\
\end{tabular}
.
\end{gather*}
\caption{Maximal PF root (the leading real root) for KL polynomials. The right plot illustrate the approximate $a\in\R$ assuming that the function on the left grows like $\sank^{a}$.}
\label{sec8-maxPF}
\end{figure}

\begin{figure}
\fcolorbox{tomato!50}{orchid!10}{\includegraphics[height=5.2cm]{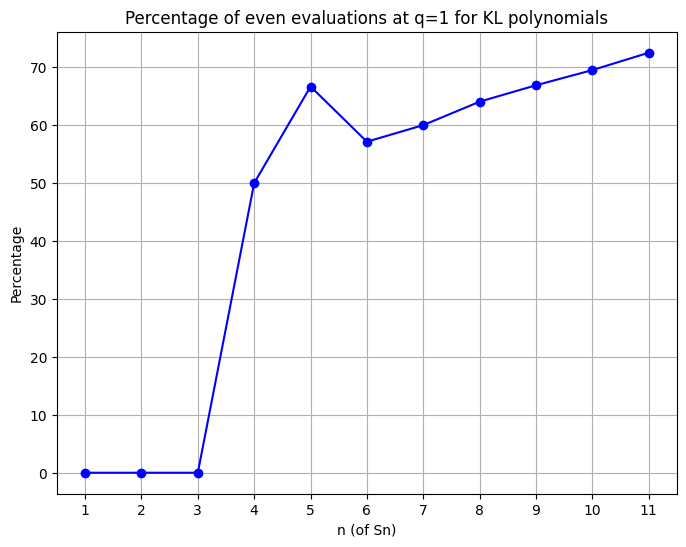}\includegraphics[height=5.2cm]{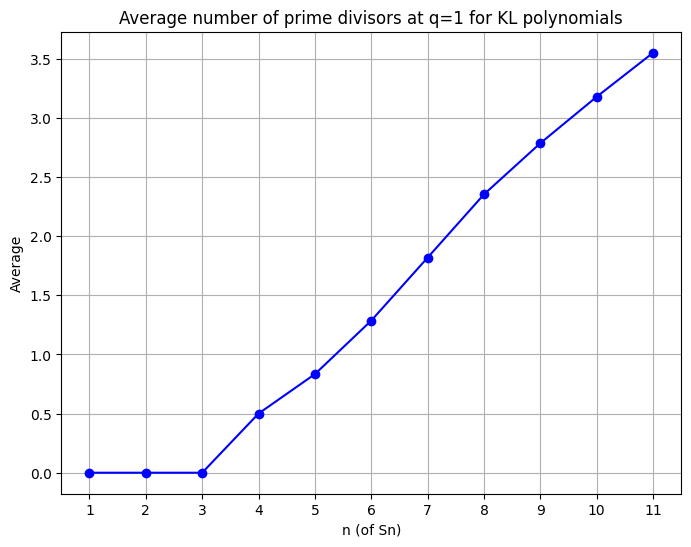}}
\begin{gather*}
\begin{tabular}{c|ccccccccccc}
$\sank$  & 1 & 2 & 3 & 4 & 5 & 6 & 7 & 8 & 9 & 10 & 11 \\
\hline
\% even & 0 & 0 & 0 & 50 & 66.67 & 57.14 & 60.00 & 64.05 & 66.88 & 69.50 & 72.50 \\
\hline
av & 0 & 0 & 0 & 0.5 & 0.8333 & 1.2857 & 1.8182 & 2.3529 & 2.7852 & 3.1777 & 3.5498 \\
\end{tabular}
.
\end{gather*}
\caption{The percentage of KL polynomials with evaluation at $\vpar=1$ being even, and the average number of prime divisors of that evaluation.}
\label{sec-Divisors}
\end{figure}

\begin{figure}[H]
\fcolorbox{tomato!50}{orchid!10}{\includegraphics[height=6.0cm]{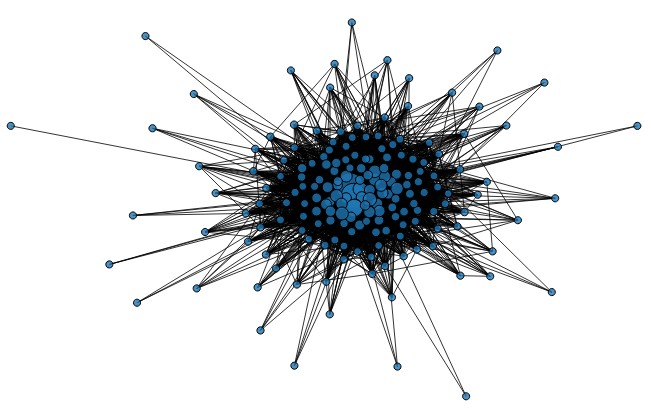}}
\caption{Ballmapper for random polynomials.}
\label{sec9-random}
\end{figure}

\begin{figure}[H]
\fcolorbox{tomato!50}{orchid!10}{\includegraphics[height=5.5cm]{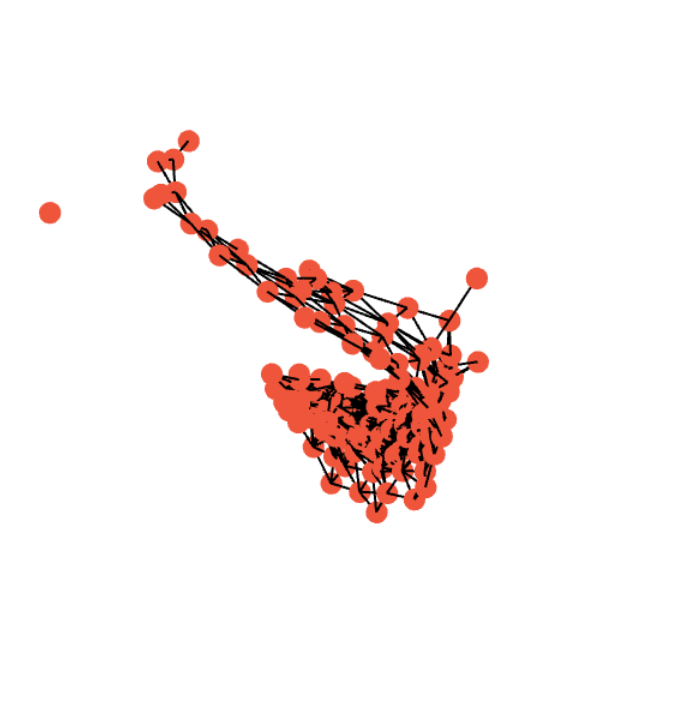}\includegraphics[height=5.5cm]{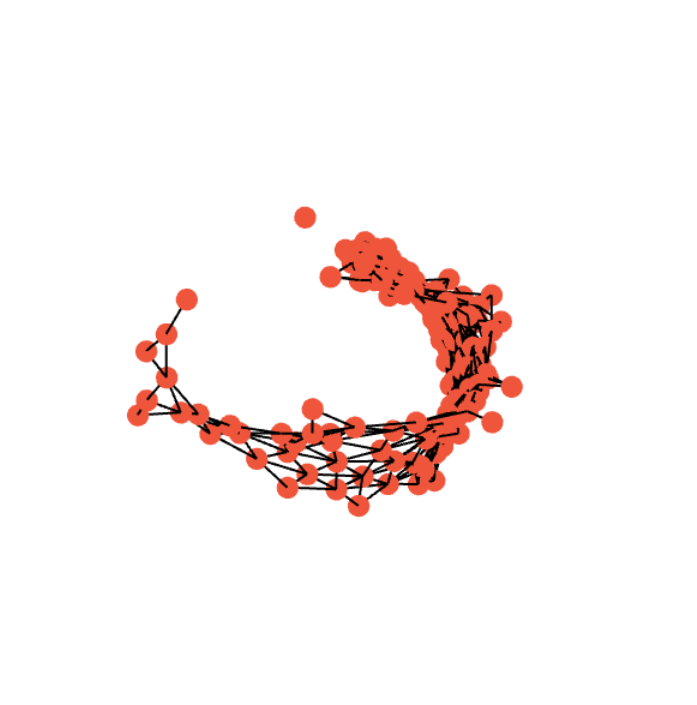}}
\fcolorbox{tomato!50}{orchid!10}{\includegraphics[height=5.5cm]{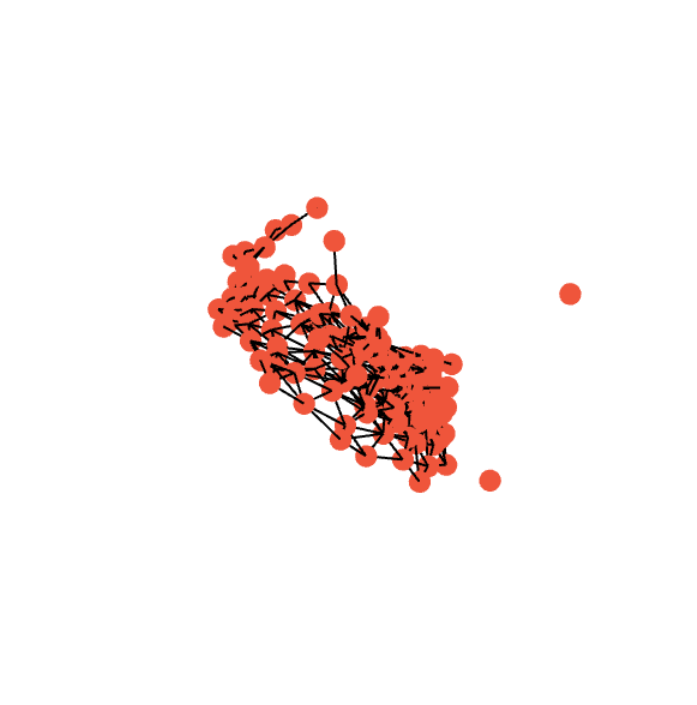}\includegraphics[height=5.5cm]{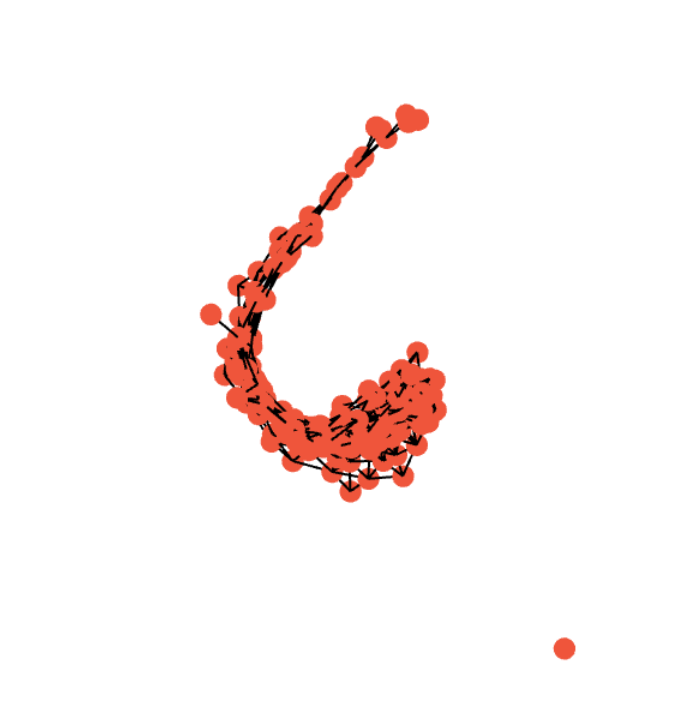}}
\caption{Spring embedding for the ballmapper graph of the $\sym[11]$ KL polynomials.}
\label{sec9-KLspring}
\end{figure}

\begin{figure}[H]
\fcolorbox{tomato!50}{orchid!10}{\includegraphics[height=5.5cm]{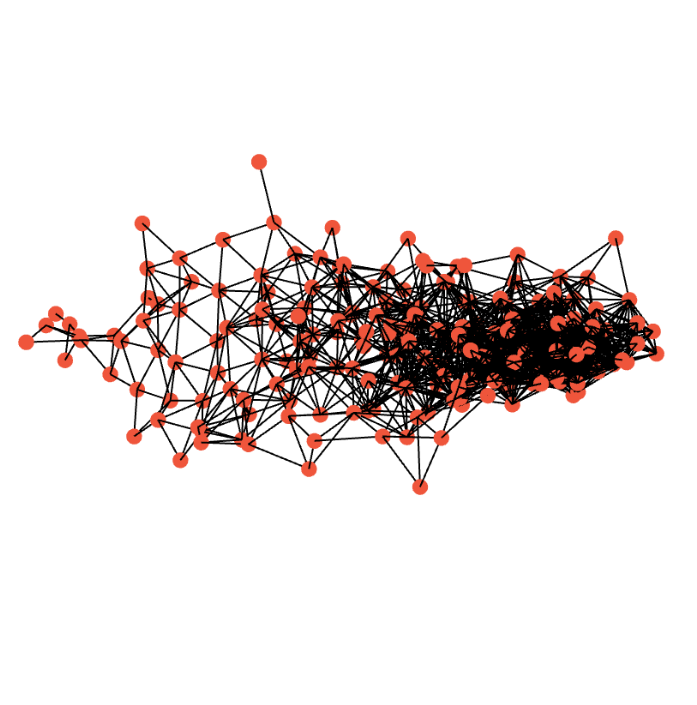}\includegraphics[height=5.5cm]{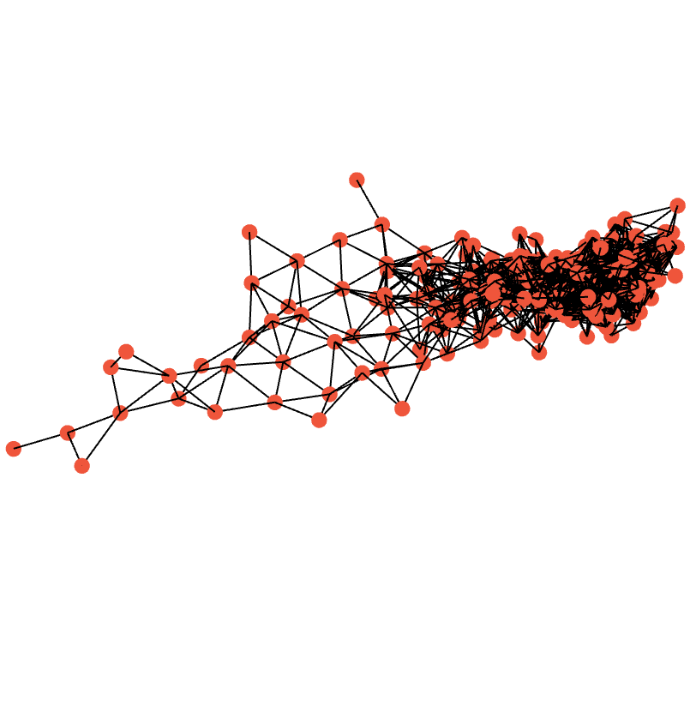}}
\fcolorbox{tomato!50}{orchid!10}{\includegraphics[height=5.5cm]{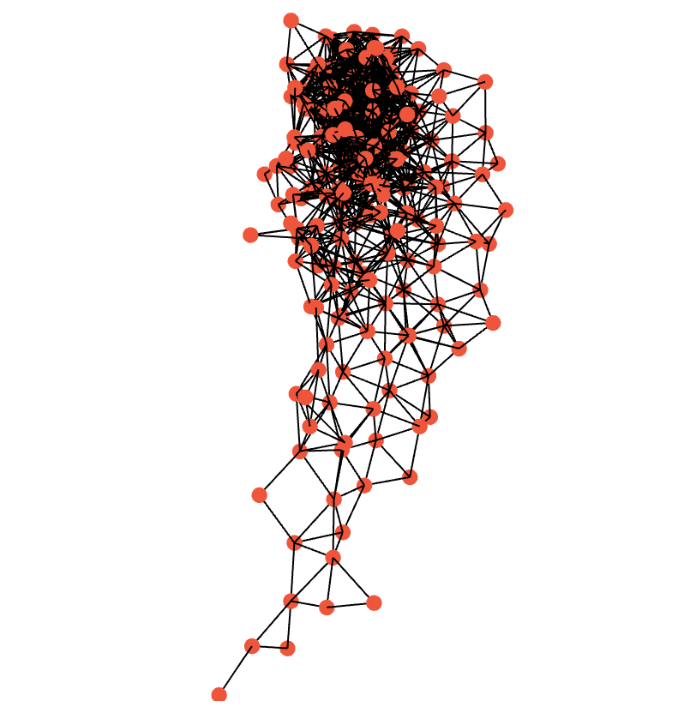}\includegraphics[height=5.5cm]{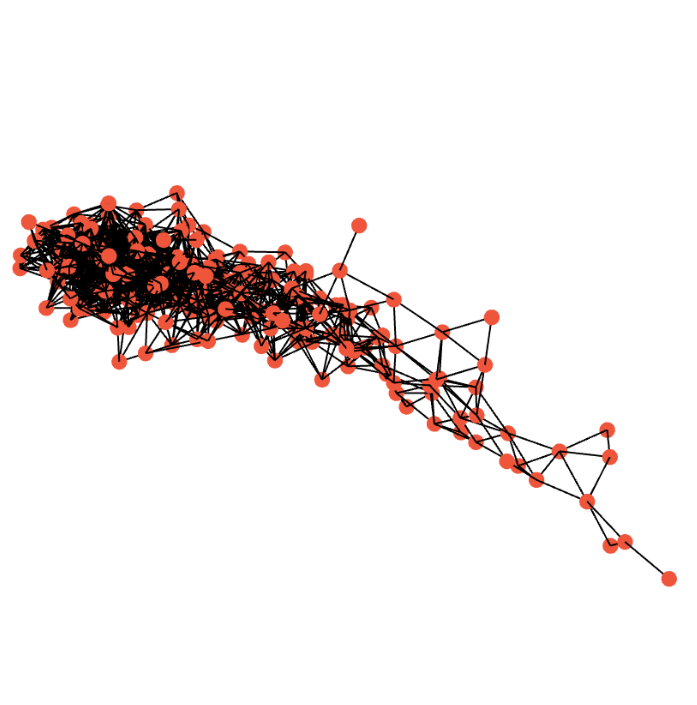}}
\caption{Kamada--Kawai embedding of the ballmapper graph of the $\sym[11]$ KL polynomials.}
\label{sec9-KLKK}
\end{figure}

\begin{figure}[H]
\fcolorbox{tomato!50}{orchid!10}{\includegraphics[height=5.5cm]{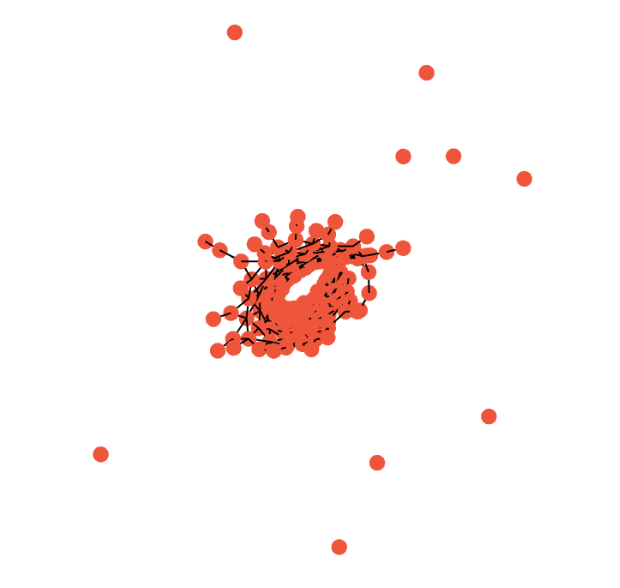}\includegraphics[height=5.5cm]{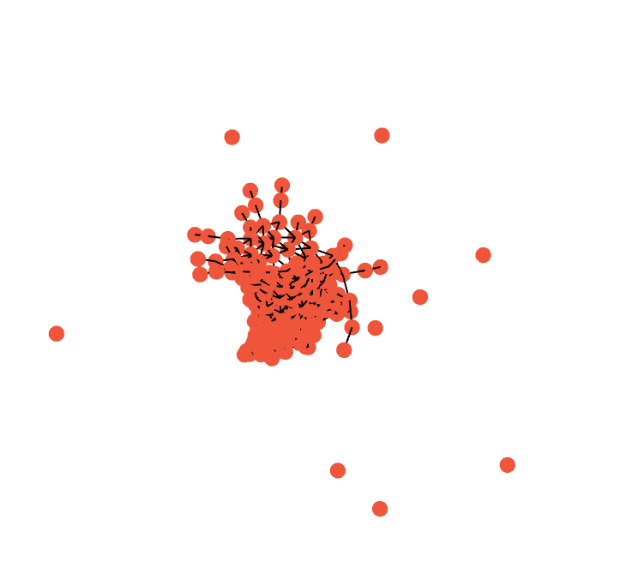}}
\fcolorbox{tomato!50}{orchid!10}{\includegraphics[height=5.5cm]{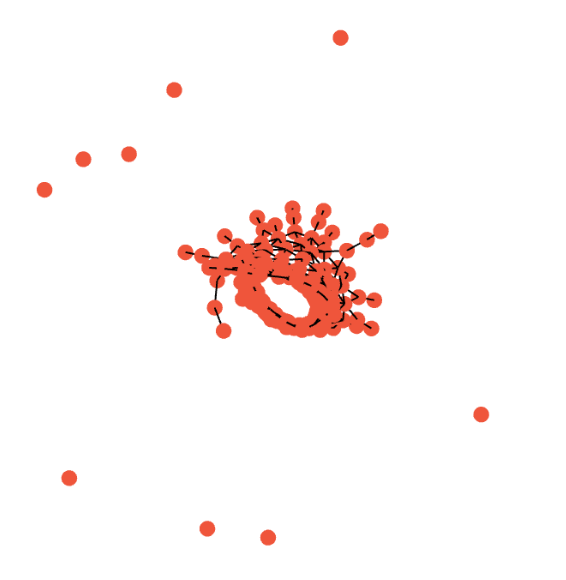}\includegraphics[height=5.5cm]{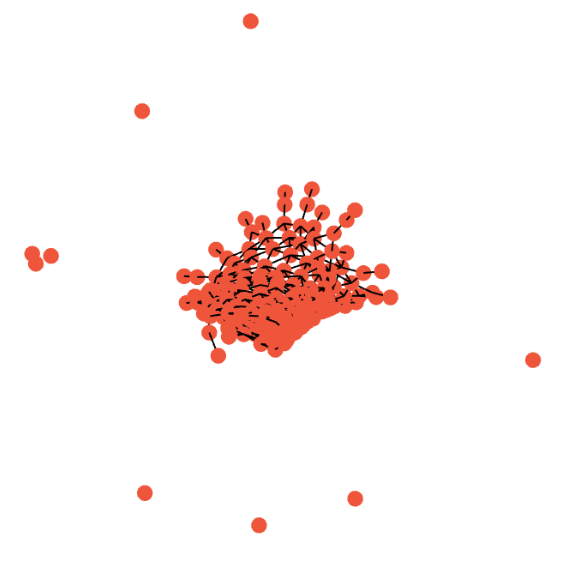}}
\fcolorbox{tomato!50}{orchid!10}{\includegraphics[height=5.5cm]{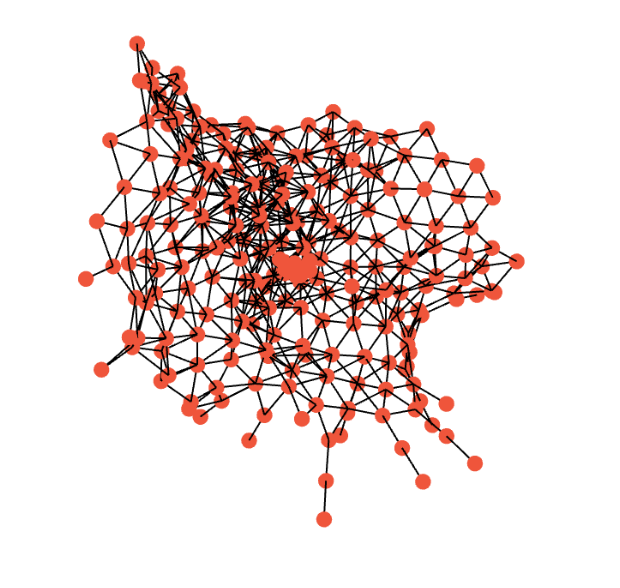}\includegraphics[height=5.5cm]{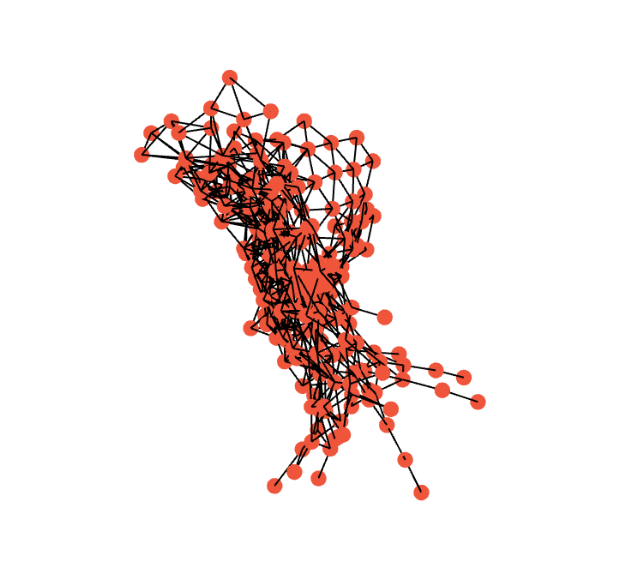}}
\fcolorbox{tomato!50}{orchid!10}{\includegraphics[height=5.5cm]{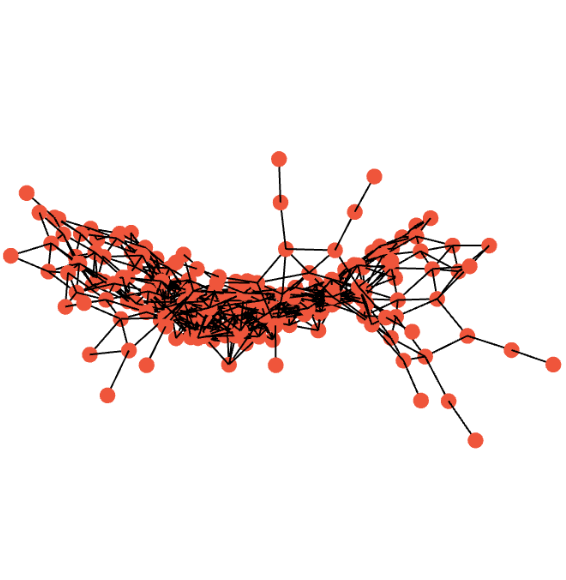}\includegraphics[height=5.5cm]{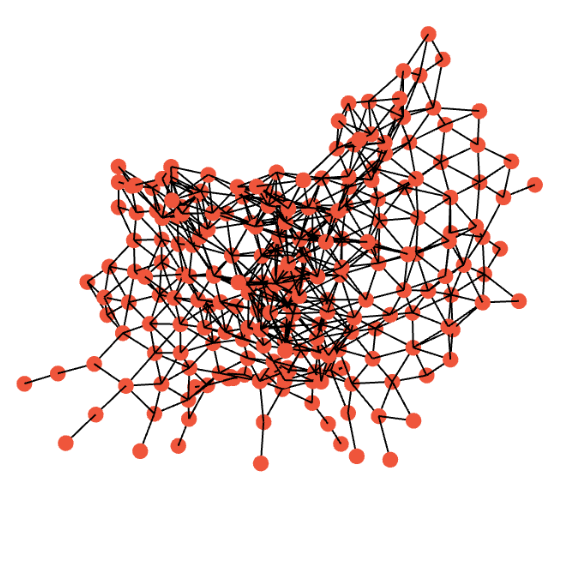}}
\caption{Spring and Kamada--Kawai embeddings for the ballmapper graph of the so-called $\sym[11]$ \emph{H polynomials} (these are essentially equivalent to the KL polynomials).}
\label{sec9-H}
\end{figure}

\newpage

\section{Preliminaries in a nutshell and notation}\label{S:Notation}

Fix the Coxeter presentation of the symmetric group $\sym$ given by simple transpositions as Coxeter generators. That is, for $m(i,i)=1$, $m(i,j)=3$ if $|i-j|=1$ and $m(i,j)=2$ if $|i-j|>1$:
\begin{gather*}
\sym=\langle
s_{1},\dots,s_{\sank-1}|(s_{i}s_{j})^{m(i,j)}=1
\rangle,
\quad
s_{i}\leftrightsquigarrow
(i,i+1).
\end{gather*}
The unit is denoted $1$.
We will illustrate elements of $\sym$ either using \emph{one-line notation} or \emph{string diagrams} (read from bottom to top), with the prototypical example being:
\begin{gather*}
\changed{\drawstringdiagramlabel{{1,5,4,7,8,9,2,10,3,6}}{10}}
\leftrightsquigarrow
(1,5,4,7,8,9,2,10,3,6)
.
\end{gather*}
\changed{This should be read, following the strings bottom to top, as ``one is mapped to one, two is mapped to five, {\etc}'' and we omit the labels from the illustrations below.}
The \emph{length} $\ell(w)$ of $w\in\sym$ is then defined using the
above Coxeter presentation \changed{as the minimal number of generators needed to write $w$}, and the reduced expressions\changed{, words with $\ell(w)$ many symbols representing $w$,} are also defined with respect to \changed{the Coxeter presentation}. Two elements $v,w\in\sym$ are in \emph{Bruhat order} with $v$ being smaller and $w$ being larger, denoted by $v\br w$, if some substring of some reduced word for $w$ is a reduced word for $v$. The Bruhat order is a partial order on $\sym$ and we will use expressions like $v\rbb w$ with the evident meaning. 
Finally, \emph{left and right descent sets} of $w\in\sym$ are
$\ldes(w)=\{s_{i}|s_{i}w\bbr w\}$ and $\rdes(w)=\{s_{i}|ws_{i}\bbr w\}$, which are used to define \emph{extremal pairs} as $\ep=\{v\br w\in\sym\times\sym|\ldes(w)\subset\ldes(v),\rdes(w)\subset\rdes(v)\}$. 
We also say $v\br w$ form an extremal pair if they are in $\ep$.

The Bruhat order can be refined into a total order, called \emph{one-line order}, with $v$ being smaller and $w$ being bigger, denoted by $v\ol w$, if $v$ is smaller than $w$ in one-line notation read lexicographically. We also use {\eg} $v\loo w$ as in the Bruhat order.

Let $\vpar$ denote some formal variable.
\cite{KaLu-reps-coxeter-groups} associates a polynomial in $\vpar$, the \emph{KL polynomial},
with pairs of elements of $\sym$. That is, one can show that we have
\begin{gather*}
v,w\in\sym
\mapsto
\KL{v}{w}=\KL{v}{w}(\vpar)\in (1+\vpar\N[\vpar])\cup\{0\}.
\end{gather*}
\ochanged{We will not repeat the
definition of the KL polynomials here but rather list some properties:}
\begin{enumerate}

\item $\KL{w}{w}=1$, and $\KL{v}{w}\neq 0$ if and only if $v\br w$.

\item For $v\bbr w$ we have $\deg\KL{v}{w}\leq\frac{1}{2}(\ell(w)-\ell(v)-1)$, and the \changed{coefficient of $\vpar^{(\ell(w)-\ell(v)-1)/2}$} is called the \emph{$\mu$ coefficient} $\mu(v,w)$. This coefficient is nonzero only if $v\br w$ form an extremal pair.

\item We have $\KL{u}{w}-\KL{v}{w}\in\N[\vpar]$ for all $u,v,w\in\sym$ with $u\br v\br w$.
In particular, $\KL{v}{w}(1)\leq\KL{1}{w}(1)$ for all $v,w\in\sym$.

\item \changed{An upper bound for t}he maximum degree is $\max\{\deg\KL{v}{w}|v,w\in\sym\}\leq\lfloor\tfrac{1}{2}\ell(w_{0})\rfloor=\lfloor\tfrac{1}{4}(\sank-1)\sank\rfloor$.

\item The computer program of our choice to compute these 
polynomials is the one from \cite{Wa-mu-for-sn}. This program is also available on \cite{LaTuVa-kl-big-data-code}. Other good programs to compute KL polynomials (which are slower in our experience but can compute them for other Coxeter groups as well) are du Cloux's program \cite{duCl-positivity-finite-hecke} and Gibson's \ochanged{(and Williamson's)} Magma programs \cite{Gi-ihecke-code}.

\end{enumerate}
All of these facts can be found in the usual literature on KL polynomials, 
{\eg} \cite{KaLu-reps-coxeter-groups,BjBr-coxeter,Lu-hecke-book}.

\begin{Remark}
Our data set seems too small to say anything beyond \cite{Wa-mu-for-sn} about the $\mu$ coefficients, but we still mentioned them via (b), as they have attracted a lot of interest in the past.
\end{Remark}

\begin{Example}\label{E:Sym3}
The following example fixes our conventions. Let $\sank=3$, and 
\begin{gather*}
\drawstringdiagram{{1,2,3}}{3}
\leftrightsquigarrow
(1,2,3)
,\quad
\drawstringdiagram{{1,3,2}}{3}
\leftrightsquigarrow
(1,3,2)
,\quad
\drawstringdiagram{{2,1,3}}{3}
\leftrightsquigarrow
(2,1,3)
,\\
\drawstringdiagram{{2,3,1}}{3}
\leftrightsquigarrow
(2,3,1)
,\quad
\drawstringdiagram{{3,1,2}}{3}
\leftrightsquigarrow
(3,1,2)
,\quad
\drawstringdiagram{{3,2,1}}{3}
\leftrightsquigarrow
(3,2,1)
.
\end{gather*}
The one-line ordering on these is $(1,2,3)\ol(1,3,2)\ol(2,1,3)\ol(2,3,1)\ol(3,1,2)\ol(3,2,1)$. In this ordering we can put the KL polynomials $\KL{v}{w}$ into a matrix:
\begin{gather*}
\begin{pmatrix}
1 & 1 & 1 & 1 & 1 & 1
\\
0 & 1 & 0 & 1 & 1 & 1
\\
0 & 0 & 1 & 1 & 1 & 1
\\
0 & 0 & 0 & 1 & 0 & 1
\\
0 & 0 & 0 & 0 & 1 & 1
\\
0 & 0 & 0 & 0 & 0 & 1
\end{pmatrix}
.
\end{gather*}
Here $v$ indexes the rows and $w$ the columns. The first row is 
the case of $\kl=\KL{1}{w}$, which we will focus on in some of the sections below.
\end{Example}

\begin{Notation}
As another piece of notation, we collect the KL polynomials as 
\emph{lists} (ordered, repetitions allowed, notation $[\placeholder]$), \emph{multisets} (not ordered, repetitions allowed, notation $(\placeholder)$) or \emph{sets} (not ordered, no repetitions, notation $\{\placeholder\}$):
\begin{gather*}
\text{list}\colon[\KL{v}{w}|v,w\in\sym]
,\quad
\text{multiset}\colon(\KL{v}{w}|v,w\in\sym)
,\quad
\text{set}\colon\{\KL{v}{w}|v,w\in\sym\}
.
\end{gather*}
For a multiset we use superscripts to encode multiplicities.
We then use the standard computer science notation to work with these, {\eg} $\sum$ for the sum over all entries, $\max$ for the maximal entry {\etc}
\end{Notation}

\begin{Example}
In \autoref{E:Sym3} we have
\begin{gather*}
(\KL{v}{w}|v,w\in\sym[3])=(0^{17},1^{19})
,\quad
\{\KL{v}{w}|v,w\in\sym[3]\}=\{0,1\}
,
\end{gather*}
and the list of KL polynomials is the given flattened matrix in that example. We have
\begin{gather*}
\#(\KL{v}{w}|v,w\in\sym[3],\KL{v}{w}\neq 0)
/\#(\KL{v}{w}|v,w\in\sym[3])=\tfrac{19}{36}
,\quad
\#\{\KL{v}{w}|v,w\in\sym[3],\KL{v}{w}\neq 0\}
/\#\{\KL{v}{w}|v,w\in\sym[3]\}=\tfrac{1}{2}
,
\end{gather*}
where, as throughout, we use \# for the number of elements, counted with \changed{multiplicity. Moreover, we have $\sum(\KL{v}{w}|v,w\in\sym[3],\KL{v}{w}\neq 0)=19$ and also $\sum(\KL{v}{w}|v,w\in\sym[3])=19$.}
\end{Example}

We now briefly explain the setting we use throughout. If not stated otherwise, we will identify a KL polynomial with a vector as follows. We use some upper bound $d$ for the degree of polynomials appearing, {\eg} $\sank^{2}$. Then, for $\KL{v}{w}\neq 0$:
\begin{gather*}
\KL{v}{w}=1+a_{1}\vpar^{1}+\dots+a_{\mcoeff}\vpar^{\mcoeff}
\leftrightsquigarrow
[1,a_{1},\dots,a_{\mcoeff}]\in\N^{\mcoeff+1}\subset\R^{\mcoeff+1}.
\end{gather*}
For $\KL{v}{w}=0$, we assign the zero vector. In this way, all KL polynomials are in the same vector space. In the language of computer science, we create the list of coefficients
$[1,a_{1},\dots,a_{\mcoeff}]$ padded with zeros until we get the desired length $\mcoeff$. Moreover, \changed{\emph{spread}} is $\sspan\KL{v}{w}=\deg \KL{v}{w}+1$.

\begin{Example}
Note that our vectors are padded with zeros, so the \changed{spread} is not necessarily equal to \ochanged{$\mcoeff+1$}. Explicitly, for $\KL{v}{w}=1+3\vpar^{5}+4\vpar^{9}$ we have $\sspan\KL{v}{w}=10$, independent of $\mcoeff$. \changed{(Note that $\mcoeff$ is chosen such that all polynomials fit into $\R^{\mcoeff+1}$; in particular, it is often much larger than the degree or the spread.)}
\end{Example}

\begin{Example}
Using the vector notation, we have $\sum\KL{v}{w}=\KL{v}{w}(1)$.
\end{Example}

\begin{Example}
For $\sym[10]$, we can choose $\mcoeff=22$\changed{, since the maximal expected degree is $\lfloor\frac{1}{4}90\rfloor=22$}. There is a KL polynomial of the form
$\KL{v}{w}=1+13\vpar+78\vpar^{2}+282\vpar^{3}+666\vpar^{4}+1068\vpar^{5}+1176\vpar^{6}+879\vpar^{7}+432\vpar^{8}+132\vpar^{9}+22\vpar^{10}+\vpar^{11}$. We get 
\begin{gather*}
\KL{v}{w}=[1,13,78,282,666,1068,1176,879,432,132,22,1,0,0,0,0,0,0,0,0,0,0,0]
\end{gather*}
as our coefficient list. We have $\sspan\KL{v}{w}=12$, $\ssum\KL{v}{w}=4750$, and $\max\KL{v}{w}=1176$.
\end{Example}

\begin{Remark}
Recall that we focus on the symmetric group only; the construction of \cite{KaLu-reps-coxeter-groups} in fact works for any Coxeter group. It would be interesting to redo the big data analysis for other Coxeter groups.
\end{Remark}

Throughout, we are interested in large $\sank$ behavior. We now fix some notation that we will use.

\begin{Notation}
For functions $f,g\colon\N\to\R_{>0}$ we use:
\begin{gather*}
\begin{aligned}
f\sim g
&\;\Leftrightarrow\;
\forall\varepsilon>0,\,\exists n_{0}
\;\text{ such that }\;
\fbox{$|\tfrac{f(n)}{g(n)}-1|<\varepsilon$},
\;\forall n>n_{0}
,
\\
f\in\Omega(g)
&\;\Leftrightarrow\;
\exists C>0,\,\exists n_{0}
\;\text{ such that }\;
\fbox{$|f(n)|\geq C\cdot g(n)$},
\;\forall n>n_{0}
,
\\
f\in O(g)
&\;\Leftrightarrow\;
\exists C>0,\,\exists n_{0}
\;\text{ such that }\;
\fbox{$|f(n)|\leq C\cdot g(n)$},
\;\forall n>n_{0}
.
\end{aligned}
\end{gather*}
We use these, in order, as \emph{asymptotic behavior, lower and upper bounds}.
\end{Notation}

\begin{Notation}
To analyze functions that behave like $\sank^{a}$ for some $a\in\R$, we use \emph{effective exponent}. For a function $f\colon\N\to\R_{>0}$, this is the function
$n\mapsto\frac{\ln(f(\sank))-\ln(f(\sank+1))}{\ln(\sank)-\ln(\sank+1)}$.
\ochanged{Note that this is $a$ for $f(\sank)=\sank^{a}$.}
\changed{(To the best of our knowledge, this is a standard concept that does not have a widely accepted name; we refer to it as the effective exponent for the reasons provided in the previous sentence.)}

For functions that behave like $a^{\sank}$ for some $a\in\R_{>1}$ we use
\emph{successive quotients}. For a function $f\colon\N\to\R_{>0}$, this is the function
$n\mapsto \frac{f(\sank+1)}{f(\sank)}$. 
\ochanged{Note that this is $a$ for $f(\sank)=a^{\sank}$.}
\end{Notation}

\begin{Notation}
We sometimes use \emph{log plots}, which means that we put a logarithmic scale on the y-axis.
\end{Notation}

\begin{Notation}\label{N:Fire}
\ochanged{We will also occasionally use \emph{heat plots} with a fiery color scheme, see {\eg} \autoref{sec3-2}. These are graphical representations of data where values are depicted within a matrix using a color gradient, black representing the coldest values (or zero in the binary version), white the hottest (or one), and a spectrum from red to yellow in between, visually emphasizing patterns, correlations, or variations across two dimensions.}
\end{Notation}

\section{\changed{Density A:} Percentage of nonzero KL polynomials}\label{S:Density}

\changed{We begin with a fundamental and somewhat predictable statistic about KL polynomials, which nonetheless provides motivation for the more unexpected results that follow.}
Recall that $\KL{v}{w}\neq 0$ if and only if $v\br w$. So, two questions that one might ask are the following:

\begin{enumerate}[label=(\thesection.\roman*)]

\item What is the asymptotic behavior $n\to\infty$ of
\begin{gather*}
\mathrm{den}_{\sank}=
\#(\KL{v}{w}|v,w\in\sym,\KL{v}{w}\neq 0)
/
\#(\KL{v}{w}|v,w\in\sym)?
\end{gather*}
(This is the density.)
That is, what is the asymptotic for the \emph{percentage} of nonzero KL polynomials. 

\item What is the asymptotic behavior $n\to\infty$ of
\begin{gather*}
\mathrm{den}^{av}_{\sank}=
\ssum(\ssum\KL{v}{w}|v,w\in\sym)/\#(\KL{v}{w}|v,w\in\sym)?
\end{gather*}
That is a similar question as before, but now we ask about the \emph{average} values of the KL polynomials at $\vpar=1$.

\item To calculate how `fiery' the pictures are \ochanged{({\cf} \autoref{N:Fire})}, we use the so-called \emph{mean (pixel) intensity} $\mathrm{fire}_{\sank}$: covert the picture to gray scale (by using $n\mapsto\frac{n-\min}{\max-\min}$) and compute the average.
What is
\begin{gather*}
\lim_{\sank\to\infty}\mathrm{fire}_{\sank}?
\end{gather*}
\end{enumerate}

\changed{We use the lexicographical order on the one-line notation to totally order permutations. We can then form two matrices: a) The binary plot matrix $B=(b_{v,w})_{v,w\in\sym}$, where $b_{vw}=1$ if $\KL{v}{w}\neq 0$, and $b_{vw}=0$ otherwise. This matrix has only entries in $\{0,1\}$ with $0$ colored black. b) We also use the heat plot matrix $H=\big(\KL{v}{w}(\vpar=1)\big)_{v,w\in\sym}$, with entries in $\N$. Using a heat map, for $\sank=2$ the picture \autoref{sec3-first} arises where the two matrices are actually the same.
The figures \autoref{sec3-1} and \autoref{sec3-2} show a few other ranks ($\sank=3$ should be compared with \autoref{E:Sym3}), with the convention that $B$ is illustrated on the left and $H$ on the right, and if all entries are zero or one, then there is only one plot, since the matrices $B$ and $H$ agree.}
\ochanged{For these, as before, the higher resolution images in \cite{LaTuVa-kl-big-data-code} may provide a clearer view of the patterns. The percentage of nonzero entries, the average entry, and the intensity are then in \autoref{sec3-3}.}

The following conjecture comes from additional data that one gets using (a) 
in \autoref{S:Notation} as the Bruhat order is much easier to compute in large scale than KL polynomials. 
\changed{In fact, by (a) 
in \autoref{S:Notation}, the conjecture is about the Bruhat order itself and one can use this to run the calculation much further than $\sank=7$ as in the plots. (See \autoref{R:GrowthBruhat} below.)
One can also start to see this in the effective exponent curve above, which starts to slow down in its decay; the curve flattens out for larger $\sank$.} For a more refined version of \autoref{O:nonzero}
we would need more data as we cannot now avoid calculating KL polynomials.

\begin{Conjecture}\label{C:nonzero}
\ochanged{For some scalar $s\in\R_{>0}$ we}
have the asymptotic behavior
\begin{gather*}
\mathrm{den}_{\sank}\sim\ochanged{s\cdot\sank^{-(2+a)}\text{ for some }a\in[0,1].}
\end{gather*}
\end{Conjecture}

We can prove the following.

\begin{Theorem}\label{C:Theorem}
We have
\begin{gather*}
\Omega\big((0.708...-\epsilon)^{\sank}\big)\ni\mathrm{den}_{\sank}\in O(\sank^{-2})
\text{ for all }\epsilon\in\ochanged{(0,0.708...)}.
\end{gather*}
Here $0.708...=\sqrt[11]{25497938851324213335/22!}$.
\end{Theorem}

\begin{proof}
\changed{The reference \cite[Theorem 1.1]{HaPi-compare-bruhat} provides the bounds as in the theorem for the elements in Bruhat order. Here the precise value for $0.708...$ can be found in \cite[Section 3]{HaPi-compare-bruhat}. Now use (a) in \autoref{S:Notation}.}
\end{proof}

In summary, some KL polynomials are much larger than the rest. \changed{The pictures in \autoref{sec3-4}, showing only the largest KL polynomials, illustrate this.}
Basically, there is nothing that one can see.
\changed{The reader might want to look at the higher resolution pictures available in \cite{LaTuVa-kl-big-data-code} to spot the very rare red and white pixels.}

This motivates:

\begin{Speculation}\label{O:nonzero}
We have
\begin{gather*}
\ochanged{\mathrm{den}^{av}_{\sank}\in\Omega(n^{-a})\text{ for \emph{all} }a\in\R_{>0}}
,\quad
\lim_{\sank\to\infty}\mathrm{fire}_{\sank}=0.
\end{gather*}
(The fire goes out, but some KL polynomials get very large.)
\end{Speculation}

For a potential way to attack \autoref{O:nonzero}, see the end of \autoref{S:Growth}.

\begin{Remark}\label{R:GrowthBruhat}
There are efficient algorithms to check \autoref{C:nonzero} because 
it is equivalent to the study of the proportion of elements $u,w\in\sym$ that are in Bruhat order; see \cite{HaPi-compare-bruhat} for a discussion of the latter. 
\changed{This was done for $\sank$ up to $110$.}
The weighted version \autoref{O:nonzero} of \autoref{C:nonzero} would need more information on
how large KL polynomials are. There are a few references on this, for example, 
\cite{BrMa-moment-graph} (which proves \autoref{S:Notation}.(c)) and 
\cite{BiBr-lower}, but not much \ochanged{beyond these references} seems to be known.
\end{Remark}

\section{\changed{Density B:} Number of KL polynomials}\label{S:Number}

From now on, we only look at $\kl=\KL{1}{w}=[1,a_{1},\dots,a_{\mcoeff}]\in\N^{\mcoeff+1}$ for $w\in\sym$, which corresponds to the first row in the pictures of \autoref{S:Density}. 
\changed{This first row exhibits unexpected behavior, as we will now explore.}

\begin{Remark}
Note that we lose almost all 
larger KL polynomials when focusing on the first row; see the final plot in \autoref{S:Density}. However, the maximal value in this final plot is $44$ and appears in the first row with the second entry being the permutation
\begin{gather*}
\drawstringdiagram{{6,7,3,4,5,1,2}}{7}
,
\end{gather*}
which \changed{has column} index $4265$. The same value appears in two other rows and the same column. This is true in general: $\KL{1}{w}$ evaluated at $\vpar=1$ is always larger than $\KL{v}{w}$ evaluated at $\vpar=1$, for all $v,w\in\sym$, see (c) in \autoref{S:Notation}.
\end{Remark}

We stress that there is a huge difference between the sets of KL polynomials and the multisets. So, the first questions we ask are:

\begin{enumerate}[label=(\thesection.\roman*)]

\item What is the asymptotic behavior $\sank\to\infty$ of
\begin{gather*}
\mathrm{num}_{\sank}=\#\{\kl|w\in\sym\}?
\end{gather*}
In other words, we ask about the asymptotic number of different KL polynomials.

\item What is the asymptotic behavior $\sank\to\infty$ of
\begin{gather*}
\mathrm{num}_{\sank}^{\%}=\#\{\kl|w\in\sym\}/
\#(\kl|w\in\sym)?
\end{gather*}
This asks for the percentage of different KL polynomials.

\end{enumerate}
\changed{The data is summarized in \autoref{sec4}. The following conjecture and its related speculation, as presented in \autoref{C:number} and \autoref{C:number2}, stem from the data. We find these results rather unexpected and lack an immediate theoretical explanation for their validity.}

\begin{Conjecture}\label{C:number}
We have \emph{superexponential growth}, {\ie}:
\begin{gather*}
\mathrm{num}_{\sank}\in\Omega(\gamma^{\sank})\text{ for \emph{all} }\gamma\in\R_{>1}.
\end{gather*}
\end{Conjecture}

\begin{Speculation}\label{C:number2}
\ochanged{We have \emph{subexponential decay}, {\ie}:}
\begin{gather*}
\ochanged{\lim_{\sank\to\infty}\sqrt[\sank]{\mathrm{num}_{\sank}^{\%}}\in(0,1].}
\end{gather*}
\end{Speculation}

\ochanged{\autoref{C:number2} is justified by the following, as well as the final few numbers in \autoref{sec4}.}

\begin{Theorem}
\ochanged{If $\sqrt[\sank]{\mathrm{num}_{\sank}}$ grows linearly up to a sublinear correction factor, then 
$\lim_{\sank\to\infty}\sqrt[\sank]{\mathrm{num}_{\sank}^{\%}}\in(0,1]$.}
\end{Theorem}

\begin{proof}
\ochanged{Note that $\#(\kl|w\in\sym)=\sank!$, which is the denominator of 
$\mathrm{num}_{\sank}^{\%}$. Stirling's formula gives $\sqrt[\sank]{\sank!}\sim\frac{1}{e}\cdot\sank$, and the claim follows.}
\end{proof}

We do not know how to prove \autoref{C:number}. In fact, 
even when we assume the combinatorial invariance conjecture (see {\eg} \cite{BlBuDaVeWi-sym} for an overview and some indication that it is true), we would still be left with understanding the Bruhat intervals of the symmetric group. 
Counting these is related to pattern avoidance, see, for example, \cite{Te-pattern-bruhat},
and appears to be a very difficult problem.

\changed{The above can be interpreted as follows: the complexity of KL polynomials is primarily concentrated in the first row, with the remaining rows essentially reflecting the first row's information distributed across statistics from the Bruhat order. With this understanding, we will now focus exclusively on the first row.}

\section{\changed{Extremes:} Growth of KL polynomials}\label{S:Growth}

\changed{When examining a distribution, it is crucial to consider its extremes, as they often reveal special and intriguing phenomena. As we will demonstrate, this seems to hold true for KL polynomials, {\cf} \autoref{C:Growth}, where some exciting geometry and topology appears. Consequently, we are interested in questions such as:}
\begin{enumerate}[label=(\thesection.\roman*)]

\item \label{Q:GrowthA} What is the asymptotic behavior $\sank\to\infty$ of
\begin{gather*}
\mathrm{ev}_{\sank}=\max\{\ssum\kl|w\in\sym\}?
\end{gather*}
That is, we ask how fast the KL polynomials evaluated at $\vpar=1$ grow in the rank $\sank$.

\item \label{Q:GrowthB} Similarly, what is the asymptotic behavior $\sank\to\infty$ of
\begin{gather*}
\mathrm{coeff}_{\sank}=\max\{\max\kl|w\in\sym\}?
\end{gather*}
That is, we ask how fast the coefficients of KL polynomials grow in the rank $\sank$.

\item We can ask the same questions as in \ref{Q:GrowthA} and \ref{Q:GrowthB} but for the average $\mathrm{ev}^{av}_{\sank}$ and $\mathrm{coeff}^{av}_{\sank}$ instead of the maximum. Here we mean 
\emph{average} in sense of
\begin{gather*}
\mathrm{ev}^{av}_{\sank}=
\ssum\{\ssum\kl|w\in\sym\}/\#\{\kl|w\in\sym\}
,\quad
\mathrm{coeff}^{av}_{\sank}=
\ssum\{\ssum\kl/\sspan\kl|w\in\sym\}/\#\{\kl|w\in\sym\}
.
\end{gather*}
Note that we do not count the padded zeros for the average, but 
only those zeros that are in the \changed{spread} of the KL polynomials.

\end{enumerate}
Here are the data that we collected, see \autoref{sec5-evaluation} and \autoref{sec5-maximum}.
\changed{The following, namely \autoref{C:Growth}, is justified by the plots, but also by the outline of proof given below. The corresponding conjecture is \autoref{A:Growth} which we checked for $\sank$ up to thirteen.}

\begin{Conjecture}\label{C:Growth}
We have \emph{superexponential growth}, {\ie}:
\begin{gather*}
\mathrm{coeff}_{\sank}\in\Omega(\gamma^{\sank})\text{ for \emph{all} }\gamma\in\R_{>1}.
\end{gather*}
Since $\mathrm{ev}_{\sank}\geq\mathrm{coeff}_{\sank}$, the same holds for 
$\mathrm{ev}_{\sank}$.
\end{Conjecture}

\begin{Speculation}\label{O:Growth}
We have \emph{at least exponential growth}, {\ie}:
\begin{gather*}
\mathrm{coeff}^{av}_{\sank}\in\Omega(\gamma^{\sank})\text{ for some }\gamma\in\R_{>1}.
\end{gather*}
Since $\mathrm{ev}^{av}_{\sank}\geq\mathrm{coeff}^{av}_{\sank}$, the same holds for $\mathrm{ev}^{av}_{\sank}$.
\end{Speculation}

We now discuss a proof strategy for \autoref{C:Growth}. For $k+2l=\sank$, we consider
\begin{gather*}
u_{l}^{k}=
\begin{tikzpicture}[anchorbase]
\def\permutation{{8, 9, 10, 4, 5, 6, 7, 1, 2, 3}} % Input permutation ({\eg}, {8, 9, 10, 4, 5, 6, 7, 1, 2, 3})
\def\n{10} % Length of the permutation
% Draw the horizontal points (top and bottom)
\foreach \i in {1,...,\n} {
\draw[fill=black] (\i, 0) circle (0pt); % Bottom row points
\draw[fill=black] (\permutation[\i-1], 2) circle (0pt); % Top row points
}
% Draw the strings (lines between points)
\foreach \i in {1,...,\n} {
% Determine the color for the first three and last three strings
\pgfmathsetmacro{\isred}{(\i<=3 || \i>=\n-2) ? 1 : 0}
\ifnum\isred=1
\draw[usual, tomato] (\i, 0) -- (\permutation[\i-1], 2);
\else
\draw[usual, black] (\i, 0) -- (\permutation[\i-1], 2);
\fi
}
% Draw underbraces with labels
\node[tomato] at (2, -0.5) {$\underbrace{\hspace{6em}}_{\textstyle l}$}; % Underbrace for first 3 strands with label "l"
\node at (5.5, -0.5) {$\underbrace{\hspace{10em}}_{\textstyle k}$}; % Underbrace for next 4 strands with label "k"
\node[tomato] at (9, -0.5) {$\underbrace{\hspace{6em}}_{\textstyle l}$}; % Underbrace for last 3 strands with label "l"
\end{tikzpicture}
,
\end{gather*}
which is the element that pulls the first $l$ strings to the last $l$ positions, fixes $k$ strings in the middle, and pulls the last $l$ strings to the first $l$ positions.
In formulas $u_{l}^{k}=(\sank-l+1,
\dots,\sank,l+1,\dots,\sank-l,1,\dots,l)$ (one-line notation). 
\changed{This is a crosshatch pair $1\leq u_{l}^{k}$ in the sense of \cite{Wa-mu-for-sn}.}

We use this element to define
\begin{gather*}
w_{l}^{k}=(w_{0}^{l}\otimes\idmor_{k}\otimes w_{0}^{l})\circ u_{l}^{k},
\end{gather*}
which takes $u_{l}^{k}$ and stacks the longest words $w_{0}^{l}$ in $\sym[l]$ on the first and last $l$ strings.

\begin{Example}
For $l=1$ we have $w_{1}^{k}=u_{1}^{k}$, for example,
\begin{gather*}
w_{1}^{8}=u_{1}^{8}=
\drawstringdiagram{{10,2,3,4,5,6,7,8,9,1}}{10}.
\end{gather*}
For $l>1$ we have $w_{1}^{k}\neq u_{1}^{k}$, for example,
\begin{gather*}
w_{2}^{6}=
\begin{gathered}
\drawstringdiagram{{2,1,3,4,5,6,7,8,10,9}}{10}
\\[-0.2cm]
\circ
\\[-0.2cm]
\drawstringdiagram{{9,10,3,4,5,6,7,8,1,2}}{10}
\end{gathered}
,
\end{gather*}
where the bottom diagram is $u_{2}^{6}$ and the top diagram is $w_{0}^{2}\otimes\idmor_{6}\otimes w_{0}^{2}$.
\end{Example}

\ochanged{Here is a conjecture that is interesting in its own right:}

\begin{Conjecture}\label{A:Growth}
We have that $\kl[w_{l}^{k}]$ is the Hilbert--Poincar{\'e} polynomial of $(\C\Pp^{l})^{\times(k-1)}$:
\begin{gather*}
\kl[w_{l}^{k}]=(1+\vpar+\dots+\vpar^{l})^{k-1}.
\end{gather*}
\end{Conjecture}

\begin{Remark}[Added after publication]\label{R:ChatGPT}
We were contacted on the 26.June.2026 by Eivind Otto Hjelle (eohjelle@gmail.com), who shared a proof of \autoref{A:Growth} generated by ChatGPT 5.5 Pro using the following prompt:

\begin{tcolorbox}[colback=green!5!white, colframe=green!10!black, title=Prompt]
Let n be an integer and k+2l=n. Using one-line notation, define a 
permutation w on n elements by w = [n,n-1,...,k+l+1,l+1,l+2,...,l+k,l,l-1,...,1]. 
That is, w is the composition of transposing the first block of size 
l with the last block of size l, and reversing the order inside each 
block. Prove that the Kazhdan-Lusztig polynomial P\_\{e, w\}(q)=(1+q+…+q\textasciicircum{}l)\textasciicircum{}\{k-1\}. 
Do not use internet.
\end{tcolorbox}

We have reviewed the generated argument and were unable to find any mistakes in the logic.

The complete proof is available on \cite{LaTuVa-kl-big-data-code}.
\end{Remark}

\begin{Theorem}\label{T:Growth}
\leavevmode
\begin{enumerate}

\item \autoref{A:Growth} implies \autoref{C:Growth}.

\item \autoref{A:Growth} is true for $l=0$ and $l=1$. This implies that 
$\mathrm{ev}_{\sank}$ and $\mathrm{coeff}_{\sank}$ grow exponentially.

\end{enumerate}
\end{Theorem}

\begin{proof}
\textit{(a).} \ochanged{Recall that the coefficients of $(x_{1}+\dots+x_{a})^{ba}$ are the a-nomial coefficients (sometimes called multinomial coefficients).} The middle a-nomial coefficients satisfy
\begin{gather*}
\binom{ba}{a,\dots,a}\sim 
\sqrt{\tfrac{6}{(a^{2}-1)\pi}}\cdot b^{-1/2}\cdot a^{b},
\end{gather*}
as follows from Stirling's approximation.
Since we can \ochanged{fix $l$ to be any $a\in\N$ we choose}, the conjecture follows.

\textit{(b).} For $l=0$ this is clear, while the $l=1$ case follows from \cite[Theorem 2]{ShShVa-kl-flags}\changed{: that theorem proves that $\kl[w_{1}^{k}]$ is the Hilbert--Poincar{\'e} polynomial of the $k$th power of $S^{2}\cong\C\Pp^{1}$ by relating the Schubert variety associated with $w_{1}^{k}$ to $(S^{2})^{\times k-1}$}. In particular, as in (a), Stirling's approximation gives 
$\mathrm{ev}_{\sank}\geq\mathrm{coeff}_{\sank}\in\Omega(\sank^{-1/2}\cdot 2^{\sank})$.
\end{proof}

Let us recall the proof of \cite[Theorem 2]{ShShVa-kl-flags}, as one might hope that it generalizes to prove \autoref{A:Growth} in general. The authors of \cite{ShShVa-kl-flags} construct a geometric space, a flag variety, and show that this space 
admits a small resolution of singularities. This, in turn, shows that the computation of KL polynomials boils down to computing ordinary cohomology, which in the $l=1$ case is $(1+\vpar)^{b}$: the cohomology of a $b$ times product of the sphere $S^{2}$, or equivalently, of $\C\Pp^{1}$.

Replacing $\C\Pp^{1}$ by $\C\Pp^{l}$, one could expect something similar to work, in general. However, in this case one does not have a small resolution of singularities, so some other argument is needed.

\begin{Example}\label{E:Big}
We first made the naive guess that $u_{l}^{k}$ can be used to generalize
\cite[Theorem 2]{ShShVa-kl-flags}. However, the KL polynomial of $u_{l}^{k}$ is rather difficult and we do not know its growth \ochanged{rate.}

In fact, the coefficients of these become really large. For example, for $\sym[13]$ and
\begin{gather*}
u_{3}^{7}=\drawstringdiagram{{11,12,13,4,5,6,7,8,9,10,1,2,3}}{13}
\end{gather*}
we have
\begin{gather*}
\kl[u_{3}^{7}]=[1,30,433,3994,26119,127617,481228,1431090,3404124,6529693,
\\
10129212,12681891,12724552,10099918,6218204,2889956,978297,229796,34889,3014,118,1].
\end{gather*}
(This polynomial took 60 days to compute). 
The largest possible coefficient from \autoref{A:Growth} below for $\sank=13$ is $1107$, indicating that \autoref{A:Growth} significantly underestimates the actual growth.
\end{Example}

\begin{Remark}\label{R:Growth2}
For this section, the same references as in \autoref{E:Big}
are relevant. Indeed, we will comment on one momentarily.
\end{Remark}

One way to address \autoref{O:nonzero} or 
\autoref{O:Growth}
could be the following. \cite[Theorem 1]{BiBr-lower} implies that, for fixed $m\in\N$, the set
\begin{gather*}
\{w\in\sym|\ochanged{{\textstyle\sum}\kl}\leq m\}
\end{gather*}
is characterized by \emph{pattern avoidance}. For example, $m=1$ corresponds to the avoidance of
\begin{gather*}
4231\colon
\drawstringdiagram{{4,2,3,1}}{4}
,\quad
3412\colon
\drawstringdiagram{{3,4,1,2}}{4}
,
\end{gather*}
\ochanged{and the $m=2$ case is }nicely explained in \cite[Appendix]{Wo-pattern}. In turn, the avoidance of these patterns is well studied and surprisingly strong results are known. For example, for the two patterns $4231$ and $3412$, we have
\begin{gather*}
\#\{w\in\sym|w\text{ avoids }4231\}
\sim
(1.15...)\cdot(0.95...)^{\sank}\cdot\sank!
,\quad
\#\{w\in\sym|w\text{ avoids }3412\}
\sim
(1.14...)\cdot(0.95...)^{\sank}\cdot\sank!
,
\end{gather*}
with precise values for the numbers, see \cite[A113228 and A113229]{Oeis}
and, for example, \cite[Section 4.2.4]{DoKh-pattern}.
Using results of this form systematically might lead to some insight regarding \autoref{O:nonzero} or
\autoref{O:Growth}.

\section{\changed{Structure A:} Unimodality of KL polynomials}\label{S:Unimodal}

\changed{To begin our analysis in \autoref{S:Unimodal} through \autoref{S:Ballmapper}, we explore the distribution of KL polynomials by comparing it to distributions of uniformly distributed polynomials, random for short, as collected in various lemmas below. Remarkably, to the best of our knowledge, the distribution of KL polynomials exhibits distinctive features that differ significantly from any known distribution. This unexpected property of KL polynomials is particularly striking, as it contrasts with the notion suggested in \cite{Po-kl-symmetric-group} that ``every polynomial is a KL polynomial,'' despite the reference not directly addressing their distribution.}

Recall that a \emph{multimodal distribution}, say with $k$ modes, is a probability distribution with $k$ local peaks of the distribution. The case $k=1$ is called \emph{unimodal}, 
$k=2$ is \emph{bimodal}, $k=3$ is called \emph{trimodal} {\etc}

Applying the same terminology to the vector $\kl=[1,a_{1},\dots,a_{\mcoeff}]\in\N^{\mcoeff+1}$ for $w\in\sym$ raises the following questions:

\begin{enumerate}[label=(\thesection.\roman*)]

\item What is the limit $n\to\infty$ of
\begin{gather*}
\mathrm{uni}_{\sank}=
\#\{\kl|w\in\sym,\kl\text{ is unimodal}\}
/\#\{\kl|w\in\sym\}?
\end{gather*}
In other words, we are interested in the percentage of unimodal
KL polynomials. We can ask the same question for the multiset: What is the limit $n\to\infty$ of
\begin{gather*}
\mathrm{uni}_{\sank}^{mult}=
\#(\kl|w\in\sym,\kl\text{ is unimodal})
/\#(\kl|w\in\sym)?
\end{gather*}

\item Similarly, for $k\in\Z_{\geq \ochanged{2}}$, let $\mathrm{multi}_{\sank}^{k}$ denote the percentage of 
$k$ modal KL polynomials that are not $j$ modal for $j<k$, taken among the KL polynomials that are not $j$ modal for $j<k$. What is
\begin{gather*}
\lim_{\sank\to\infty}\mathrm{multi}_{\sank}^{k}?
\end{gather*}
For example, 
$\mathrm{multi}_{\sank}^{2}$ is the percentage of bimodal KL polynomials among all 
nonunimodal KL polynomials. We could track the same statistics for the multiset version.

\end{enumerate}
\ochanged{For the comparison with a random distribution, let} us first point out that there are essentially no 
unimodal polynomials (similarly for bimodal {\etc}, but we did not find any more explicit data):

\begin{Lemma}\label{L:UnimodalCount}
\changed{Let $u_{\sank}=\#\{P|P\in\Z_{>0}^{\infty},\ssum P=\sank,\kl\text{ is unimodal}\}$ and $p_{\sank}=\#\{P|P\in\Z_{>0}^{\infty},\ssum P=\sank\}$. We have the limit}
\begin{gather*}
\changed{\lim_{\sank\to\infty}\sqrt[\sank]{\frac{u_{\sank}}{p_{\sank}}}=\frac{1}{2}.}
\end{gather*}
\changed{(A precise asymptotic is given in the proof.) That is, the percentage of unimodal $[a_{0},\dots,a_{\mcoeff}]\in\Z_{>0}^{\infty}$ goes to zero \emph{exponentially} fast.}
\end{Lemma}

\begin{proof}
\changed{The two sets in the statement are the number of unimodal compositions and the number of compositions of $\sank$, which we will now estimate.
The number of compositions of $\sank$ is $2^{\sank-1}$, which follows immediately by mapping them bijectively to $\sank-1$ bit binary numbers: 
the composition $a_{1}+\dots+a_{\sank}$ corresponds to the string 
$1^{a_{1}-1}01^{a_{2}-1}0\dots 1^{a_{\sank-1}-1}01^{a_{\sank}-1}$, where the exponents indicate number of copies. Combine this with \cite[Theorem 2]{Au-stacks}, which proves that the number of unimodal compositions is asymptotic to $s\cdot n^{-1}\cdot t^{\sqrt{n}}$ for $s=1/(8\sqrt{3})$ and $t=\exp(\sqrt{2/3}\pi)$. The ratio is now straightforward to analyze.}
\end{proof}

\begin{Remark}
We expect something similar to be true for bimodal {\etc}, but we did not find any more explicit statement. However, we will use a different way to generate polynomials below, but we have not found a percentage count for these. We still see \autoref{L:UnimodalCount} and the data below as an indication that the KL polynomials have a remarkable distribution.
\end{Remark}

\begin{Lemma}
\changed{(This is not intended to be a precise statement.) Assume that KL polynomials for $\sym$ form a random distribution. Then $\lim_{\sank\to\infty}\mathrm{uni}_{\sank}=\lim_{\sank\to\infty}\mathrm{uni}_{\sank}^{mult}=0$ exponentially fast. Similarly for bimodal {\etc}}
\end{Lemma}

\begin{proof}\changed{(Sketch.)
By \autoref{L:UnimodalCount}.}
\end{proof}

\changed{The data is summarized in \autoref{sec6-unimodal} and \autoref{sec6-multimodal}.}
\ochanged{The data for the multiset version} is not spelled out because the data set is too small
and the percentages are essentially all zero for $\sym[11]$.

\begin{Conjecture}\label{C:Unimodal}
\emph{Almost all} KL polynomials are unimodal, {\ie}:
\begin{gather*}
\lim_{\sank\to\infty}\mathrm{uni}_{\sank}=\lim_{\sank\to\infty}\mathrm{uni}_{\sank}^{mult}=1.
\end{gather*}
\end{Conjecture}

\begin{Remark}
Given that the fire seems to go out as in \autoref{O:nonzero}, one could expect most KL polynomials to be rather small and thus to have a higher percentage of unimodal \ochanged{polynomials among them}. In particular, one could say that $\lim_{\sank\to\infty}\mathrm{uni}_{\sank}=1$ is more surprising than its multiset version.
\end{Remark}

The following is not fully justified by the above data, as {\eg} there are more trimodal KL polynomials for $\sym[11]$ than bimodal ones. However, we think that the pattern is that they all peak shortly after they appear, which skews the picture we get for the data we have. So:

\begin{Speculation}\label{O:Unimodal}
\emph{Almost all} KL polynomials that are not unimodal are bimodal, 
\emph{almost all} KL polynomials that are not unimodal or bimodal are trimodal, {\etc}, {\ie}:
\begin{gather*}
\lim_{\sank\to\infty}\mathrm{multi}_{\sank}^{k}=1.
\end{gather*}
Similarly for the multiset versions.
\end{Speculation}

We sadly do not see any way to prove \autoref{C:Unimodal} or \autoref{O:Unimodal} at the 
time of writing this paper. However, although we do not know any connection, related effects 
are known in the literature. Most notable what 
one could call log concavity of R polynomials, see 
{\eg} \cite[Conjecture 2.4]{Br-conjectures} for details and a precise formulation, and, for example, \cite[Corollary 5.5.3]{BjBr-coxeter} for the relation to KL polynomials.

\section{\changed{Structure B:} Roots of KL polynomials: all roots}\label{S:Roots}

We now compare the following two pictures, following the ideas in \cite{BaChDe-roots} and the references to various blogs therein. For other root plots of `famous' polynomials see {\eg} \cite{WuWa-jones-roots} (albeit this reference takes a quite different perspective).

\ochanged{\autoref{sec7-KL}} is 
the multiset of roots of the set of KL polynomials. In formulas, we plot
\begin{gather*}
(\mathrm{roots}(p)|p\in\{\kl|w\in\sym[11]\})
\end{gather*}
with the brightness indicating the density (bright=high density; as before, this is a heat plot).
\ochanged{\autoref{sec7-random}} is of the same type, but with a different set of polynomials in $1+\N[\vpar]$. \ochanged{The plot was generated using} $293189$ polynomials in $1+\N[\vpar]$, with a randomly chosen degree \ochanged{in $\{1,\dots,15\}$, and then randomly chosen coefficients $\{0,\dots,61582\}$. These numbers are chosen in this way because the
total number of distinct KL polynomials for $\sym[11]$ is $293189$, the maximal degree is $15$, and the maximal coefficient is $61582$.}
In both cases, the right plot is a zoomed plot with a marked unit circle.

The plots are very different. To make this into a clearer statement, let us first recall a few facts about the distribution of
roots of random polynomials.

\begin{Lemma}
\ochanged{(This is not intended to be a precise statement.) Assume that our sample of KL polynomials for $\sym[11]$ is a random distribution. Then:}
\begin{enumerate}

\item \ochanged{The \emph{expected percentage} of real roots is $\approx 17.23\%$.}

\item \ochanged{For $\sank\to\infty$, \emph{almost all} roots will have an absolute value in $[1-\epsilon,1+\epsilon]$ for a fixed $\epsilon\in\R_{>0}$.}

\end{enumerate}
\end{Lemma}

\begin{proof}\ochanged{(Sketch.)}
\textit{(a).} By \cite{Ka-roots-random}, the number of expected real roots of a polynomial of degree $k$ is $2/\pi\ln k$. To simplify our calculation, assume that the average degree of our sample is $e^{2}$ for the usual $e\approx2.71...$. Then the expected percentage of real roots is $4\pi/e^{2}\approx 17.23\%$.

\textit{(b).} This is expected for random polynomials; see, for example, \cite{ShVa-roots-random}. Fairly explicit formulas for the distribution are known, for example, see
\cite{MeBeFoMaAa}, but the only thing we notice here is that, in the limit, almost all roots will have an absolute value in $[1-\epsilon,1+\epsilon]$ for a fixed $\epsilon\in\R_{>0}$.
\end{proof}

\begin{Remark}
The references \cite{Ka-roots-random,ShVa-roots-random} are for polynomials with complex coefficients since we have not found any references for integer-valued polynomials.
\end{Remark}

\ochanged{Notable beyond that are the following:}
\begin{enumerate}[label=(\thesection.\alph*)]

\item The tail of negative real numbers. The tail is explained in \autoref{S:PF} below.

\item The many holes in the KL plot. We do not know a general statement, but see \cite{BaChDe-roots} for similar patterns (this should be true for other integer valued polynomials as well). The case of coefficients in $\{-1,0,1\}$ is addressed, for example, in \cite{CaKoWa}.
The only exception is `forbidden region' (the black pizza slice around zero for real value $>0$): this is true in general, and appears in both pictures, and is a classical result; see {\eg} \cite[Lemma 2]{Fi-irreducibility-nonnegative}.

\end{enumerate}

\ochanged{A} few statistics for the roots of KL polynomials \ochanged{are given in \autoref{sec7-real} and \autoref{sec7-distance}.}
As one can see, these are very different from the random distribution. So we ask:

\begin{Question}\label{Q:Distribution}
What is a \emph{random} KL polynomial (defined so that one gets values closer to the above patterns)?
\end{Question}

We do not know how to attack this problem beyond studying and adjusting the above references.

\section{\changed{Structure C:} Roots of KL polynomials: Perron--Frobenius root}\label{S:PF}

We now address the tail \ochanged{mentioned in} \autoref{S:Roots}.

A polynomial $f=[a_{0},\dots,a_{\mcoeff}]\in 1+\vpar\Z[\vpar]$ satisfies the \emph{(negative) Perron--Frobenius (PF for short) property}, or simply \emph{is PF}, if there exists a root $\lambda\in\R_{\leq 0}$ of $f$ with $-\lambda>\max\{|\mu|\mid\mu\text{ is a root of }f,\mu\neq\lambda\}$. Such a $\lambda$ is called the \emph{PF root} of $f$. The name is explained by the proof of the following statement.

\begin{Lemma}\label{L:PFroot}
Almost all $f\in 1+\vpar\Z[\vpar]$ satisfy the PF property, {\ie}:
\begin{gather*}
\lim_{M\to\infty}
\frac{\#\{f\in 1+\vpar\Z[\vpar]|\max[a_{0},\dots,a_{\mcoeff}]\leq M,\mcoeff\leq M,f\text{ is PF}\}}
{\#\{f\in 1+\vpar\Z[\vpar]|\max[a_{0},\dots,a_{\mcoeff}]\leq M,\mcoeff\leq M\}}
=1.
\end{gather*}
\end{Lemma}

\begin{proof}
Let $\bar{f}$ be $f$ normalized, {\ie} if $f=[1,a_{1},\dots,a_{\mcoeff}]$, $a_{\mcoeff}\neq 0$, then 
$\bar{f}=[1/a_{\mcoeff},a_{1}/a_{\mcoeff},\dots,1]$. Denote its coefficients by $b_{i}$.
The companion matrix of $f$, having $\bar{f}$ as its characteristic polynomial, is
\begin{gather*}
\begin{psmallmatrix}
0 & 0 & 0 & \cdots & 0 & -b_{0} \\
1 & 0 & 0 & \cdots & 0 & -b_{1} \\
0 & 1 & 0 & \cdots & 0 & -b_{2} \\
0 & 0 & 1 & \cdots & 0 & -b_{3} \\
\vdots & \vdots & \vdots & \ddots & \vdots & \vdots \\
0 & 0 & 0 & \cdots & 1 & -b_{\mcoeff-1} \\
\end{psmallmatrix}.
\end{gather*}
Ignoring the sign in the final column, the associated graph is strongly connected, unless some $b_{i}=0$ (and this is an if and only if),
and the Perron--Frobenius theorem applies. Taking the sign then back into account gives the negative Perron--Frobenius eigenvalue, and therefore root. Finally, vanishing of a coefficient occurs with probability $1-(\tfrac{M-1}{M})^{M}$, but the companion matrix will still satisfy the PF property with one exception: if $b_{\mcoeff-1}=0$. As this happens with probability $\tfrac{1}{M}$, the limit $M\to\infty$ of the converse condition, as in the statement, is $1$.
\end{proof}

\ochanged{Since KL polynomials $\kl$ are in $1+\vpar\Z[\vpar]$, we can ask:}
\begin{enumerate}[label=(\thesection.\roman*)]

\item What is the limit $n\to\infty$ of
\begin{gather*}
\mathrm{PF}_{\sank}=
\#\{\kl|w\in\sym,\kl\text{ is PF}\}
/\#\{\kl|w\in\sym\}?
\end{gather*}
Similarly, but for the multiset, we would like to know the limit 
$n\to\infty$ of
\begin{gather*}
\mathrm{PF}_{\sank}^{mult}=
\#(\kl|w\in\sym,\kl\text{ is PF})
/\#(\kl|w\in\sym).
\end{gather*}

\item Also, what is the asymptotic behavior $\sank\to\infty$ of
\begin{gather*}
\mathrm{PFmax}_{\sank}=
\max\{-\lambda|\lambda\text{ is PF root of }\kl,w\in\sym,\kl\text{ is PF}\}?
\end{gather*}
(The sign ensures that this is a nonnegative real number.)

\end{enumerate}

\begin{Lemma}
\ochanged{(This is not intended to be a precise statement.) Assume that our sample of KL polynomials for $\sym[11]$ is a random distribution. Then, \emph{almost surely}, we have $\mathrm{PFmax}(\text{random})\geq 248.1572$.}
\end{Lemma}

\begin{proof}
\ochanged{(Sketch.)} For random polynomials {\eg} \cite{Ge-spectral,Ba-circular-law} implies that the PF root gets very large. In particular, almost surely we have $\mathrm{PFmax}(\text{random})\geq\sqrt{M}$ where $M$ is the bound for the coefficients. In \autoref{S:Roots} we have $M=61582$, so $\mathrm{PFmax}(\text{random})\geq 248.1572$.
\end{proof}

Moreover, \autoref{A:Growth} would imply that $\mathrm{PFmax}(\text{random})$ grows superexponentially if the KL polynomials were randomly distributed.

\ochanged{We summarize the data related to the PF property and to the maximal PF root in \autoref{sec8-PF} and \autoref{sec8-maxPF}.}
One pattern seems to be clear, so:

\begin{Conjecture}\label{C:PF}
\emph{Almost all} KL polynomials satisfy the PF property, {\ie}:
\begin{gather*}
\lim_{\sank\to\infty}\mathrm{PF}_{\sank}=1
,\quad
\lim_{\sank\to\infty}\mathrm{PF}^{mult}_{\sank}=1.
\end{gather*}
\end{Conjecture}

Let us motivate \autoref{S:GrowthPF}, which is up next. Assuming \autoref{C:Unimodal}, almost all KL polynomials are unimodal. In fact, the
ones with the largest coefficients tend to be unimodal as well. However, it seems that the roots of unimodal polynomials grow slowly. The most famous example of this
is $(1+\vpar+\dots+\vpar^{l})^{k-1}$, which we have seen in \autoref{S:Growth}: \ochanged{t}heir roots are roots of unity.
A nontrivial example
is \emph{chromatic roots}, the roots of chromatic polynomials (that these polynomials are unimodal is a nontrivial theorem of \cite{Hu-roots}). Since the complete graph with $k$ vertices has a chromatic root at $k-1$, the rate of growth of the chromatic roots is linear in the number of vertices.

\begin{Speculation}\label{S:GrowthPF}
We have \emph{subexponential growth}, {\ie}:
\begin{gather*}
\mathrm{PFmax}_{\sank}\in O(\gamma^{\sank})\text{ for \emph{all} }\gamma\in\R_{>1}.
\end{gather*}
\end{Speculation}

Being probably related to \autoref{S:Unimodal}, we do not know how to prove \autoref{C:PF} or \autoref{S:GrowthPF}.

\section{\changed{Structure D: Prime divisors of KL polynomials}}\label{S:Divisors}

\changed{Since $\sum\kl$ (the evaluation at $\vpar=1$) is an integer, we can ask:}
\begin{enumerate}[label=(\thesection.\roman*)]

\item \changed{What is the limit $n\to\infty$ of}
\begin{gather*}
\changed{\mathrm{even}_{\sank}=
\#\{\kl|w\in\sym,{\textstyle\sum}\kl\in2\Z\}
/\#\{\kl|w\in\sym\}?}
\end{gather*}
\changed{That is, for how many $w\in\sym$ is $\sum\kl$ even (and if true, we call $\kl$ itself \emph{even})?}

\item \changed{Also, what is the asymptotic behavior $\sank\to\infty$ of}
\begin{gather*}
\changed{\mathrm{pr}_{\sank}^{av}=
\#(p|p\text{ prime and divides }\ssum\kl,w\in\sym)
/\#\{\kl|w\in\sym\}}
\end{gather*}
\changed{This measures the average number of prime divisors of $\sum\kl$, counted with multiplicity.}

\end{enumerate}
\changed{Here is the comparison to random polynomials:}

\begin{Lemma}\label{L:Divisors}
\changed{(This is not intended to be a precise statement.) Assume that KL polynomials for $\sym$ form a random distribution. Then:}
\begin{enumerate}

\item \changed{We have $\lim_{n\to\infty}\mathrm{even}_{\sank}=\frac{1}{2}$.}

\item \changed{For $\sank=11$, the \emph{expected average} number of prime divisors is $\approx 2.9678$.}

\end{enumerate}
\end{Lemma}

\begin{proof}\changed{(Sketch.)
\textit{(a).} Clear.}

\changed{\textit{(b).} Using the notation and results of \cite[Appendix A]{KT}, the expected average number of prime divisors (counted with multiplicity) $\frac{1}{N}\sum_{x\leq N}\Omega(x)$ of a number $N$ is logarithmic in the number of digits. More precisely, $\frac{1}{N}\sum_{x\leq N}\Omega(x)=\ln\ln N+1.03\dots+O(1/\ln N)$, with $1.03\dots$ explicitly given in \cite[Appendix A]{KT}.
Now we just use this for the average evaluation of KL polynomials at $\vpar=1$
as in \autoref{sec5-evaluation}.}
\end{proof}

\begin{Remark}
\changed{The $2.9678$ from \autoref{L:Divisors}.(b) looks 
close to the $3.5498$ we got in \autoref{sec-Divisors}, but, as in the proof of \autoref{L:Divisors}, they are actually quite different as the expected scale is double logarithmic.}
\end{Remark}

\changed{The data is summarized in \autoref{sec-Divisors}.
Based on the data, we conjecture that, very surprisingly, almost all KL polynomials are even, and we speculate that the average number of prime divisors grows essentially linearly (not double logarithmic):}

\begin{Conjecture}\label{C:Primes}
\changed{\emph{Almost all} KL polynomials are even, {\ie}:}
\begin{gather*}
\changed{\lim_{\sank\to\infty}\mathrm{even}_{\sank}=1.}
\end{gather*}
\end{Conjecture}

\begin{Speculation}\label{C:Primes2}
\changed{We have \emph{superlogarithmic growth}, {\ie}}
\begin{gather*}
\changed{\mathrm{pr}_{\sank}^{av}\in\Omega\big(\log(\sank)^{\gamma}\big)\text{ for \emph{all} }\gamma\in\R_{\geq 1}.}
\end{gather*}
\end{Speculation}

\changed{We do not know any way to prove \autoref{C:Primes} or \autoref{C:Primes2}. However, these patterns look similar to facts about characters of symmetric group as {\eg} in 
\cite{Mi}, and one might hope for a connection.}

\section{\changed{Structure E:} The KL ballmapper}\label{S:Ballmapper}

Mapper algorithms are classical tools in TDA, originally introduced in \cite{mapper}, designed to explore and visualize data.
These algorithms combine dimensionality reduction, clustering, and graph network techniques to create a graph out of data. We will use a mild modification that is a bit more suitable for our purpose, as it gives nicer visualizations, \emph{ball mapper} from \cite{Dl-ballmapper}. These types of algorithms have found important applications, even in the real-world sciences; see, for example, 
\cite{Dl-tda} for references and examples of applications.
Moreover, the current section is inspired by \cite{DlGuSa-mapper}, which uses the ball mapper on data from knot theory, Tic-Tac-Toe games, and classical cancer databases. We now explain 
the results if one applies these methods to the KL polynomials.

In fact, following \cite{BlBuDaVeWi-sym}, we also used the so-called \emph{H polynomials} together with the KL polynomials. See \cite[Section 3.1]{BlBuDaVeWi-sym} for the definition. The point is that knowing the tuple $\big(\kl,\ell(w)\big)$ is equivalent to knowing the H polynomial for $w$ used in the ball mapper.

\begin{Remark}
The data sets for the H polynomials can also be found in \cite{LaTuVa-kl-big-data-code}.
\end{Remark}

We will not recall the definition of the mapper and ballmapper algorithms. Here is a very brief summary of what the reader needs to know; see the above references for details:

\begin{enumerate}

\item The algorithm we used creates a graph $G=G(\epsilon)$ 
from a point cloud using balls of a prefixed radius 
$\epsilon\in\R_{>0}$. All points in a ball with this radius are collapsed into a vertex of $G$, {\ie} $G$ which has more vertices the smaller $\epsilon$. The main point then is the study of the family of graphs when one varies $\epsilon$.

\item In order to improve the graph plots, we rescaled the KL and H polynomials,
followed by a principal component analysis (PCA).

\item For visualization in 3D,
we used force-directed drawing algorithms of (coil-)spring-type 
and Kamada--Kawai-type as provided by the Python library Bokeh.
Details about these graph visualization algorithms can be found, for example, in \cite[Chapter 12]{Ta-handbook}.

\item More details are explained in \cite{LaTuVa-kl-big-data-code}. We give only a sample of what can be found in \cite{LaTuVa-kl-big-data-code}; in particular, the graphs are much more impressive in the 3D plot in \cite{LaTuVa-kl-big-data-code}.

\end{enumerate}

\changed{For the following lemma see also \autoref{sec9-random}.}

\begin{Lemma}\label{L:RandomGraph}
\changed{(This is not intended to be a precise statement.) Assume that KL polynomials for $\sym$ form a random distribution. Then the ballmapper is a random graph.}
\end{Lemma}

\begin{proof}
\changed{(Sketch.) Immediate from the definitions (and easy to plot).}
\end{proof}

\ochanged{We now comment on the KL ball mappers for $\sym[11]$. 
The spring embedding for $\epsilon=4\cdot 10^{-9}$ and the KL polynomials, rotated to four different angles is presented in \autoref{sec9-KLspring}.
The Kamada--Kawai embedding for the the same $\epsilon$ and the KL polynomials is displayed in \autoref{sec9-KLKK}, and again in four different angles.}

We point out, as the reader is invited to check using the code in \cite{LaTuVa-kl-big-data-code}, that these arc-like (top) or fish (bottom) shapes are quite stable when varying $\epsilon$, but also when varying $\sank$ from roughly $\sank=8$ onward. This might indicate some stabilization process that happens for $\sank\to\infty$ but, at the same time, it underlines that the small $\sank$ does not provide a sample representative of KL polynomials.

Similarly, the H polynomials in spring and the Kamada--Kawai embedding, both for the same $\epsilon=10^{-9}$, are \ochanged{given in \autoref{sec9-H}.}
The same comments as for the KL polynomials apply, but this time for this web-type shape. 

ChatGPT \ochanged{model GPT-4} says about the KL polynomial ball mapper and the H polynomial ball mapper, both for the spring embedding: 
\begin{enumerate}[label=$\blacktriangleright$]

\item ``The overall shape resembles a `fish' or a `stream,' with the dense cluster forming a 'head' and the sparser connections on the left resembling a `tail.' It gives the impression of directional flow or gradient density across the network.''

\item ``The overall structure looks like a dense, irregular web, with some areas more interconnected than others. The shape does not conform to a specific geometric form, but appears somewhat sprawling and organic, indicative of a complex network with a dense central region and sparser peripheral connections.''

\end{enumerate}

Finally, again there is a huge difference between randomly generated polynomials as in \autoref{S:Roots} and KL polynomials\ochanged{, {\cf} \autoref{L:RandomGraph}.}
Thus, in analogy to \autoref{Q:Distribution}, we ask:

\begin{Question}\label{Q:DistributionTwo}
What is a \emph{random} KL polynomial (defined so that one gets BM closer to the above)?
\end{Question}

\begin{Remark}
We also tried one of the most classical tools from TDA on the KL polynomials: persistent homology. However, the data set of KL polynomials starts to be interesting very late and is then too large: we were not able to extract meaningful results.
\end{Remark}

\newcommand{\etalchar}[1]{$^{#1}$}

\end{document}